\documentclass[preprint,11pt]{elsarticle}
\usepackage{amsmath,amsfonts,amssymb}
\usepackage{indentfirst,latexsym,bm}
\usepackage{amsthm}  
\usepackage{mathrsfs} 
\usepackage{color,xcolor}
\usepackage[all]{xy}

\newtheorem{theorem}{Theorem}[section]  

\newtheorem{lemma}[theorem]{Lemma}
\newtheorem{remark}[theorem]{Remark}
\newtheorem{definition}[theorem]{Definition}
\newtheorem{problem}[theorem]{Problem}
\newtheorem{corollary}[theorem]{Corollary}

\theoremstyle{definition} 
\theoremstyle{remark}     
\theoremstyle{problem}    
\theoremstyle{question}   

\numberwithin{equation}{section}

\journal{XXX}

\begin{document}

\begin{frontmatter}





\title{The matched projections of idempotents on Hilbert $C^*$-modules}
\author[shnu]{Xiaoyi Tian}
\ead{tianxytian@163.com}
\author[shnu]{Qingxiang Xu }
\ead{qingxiang$\_$xu@126.com}
\author[HS]{Chunhong Fu}
\ead{fchlixue@163.com}
\address[shnu]{Department of Mathematics, Shanghai Normal University, Shanghai 200234, PR China}
\address[HS]{Health School Attached to Shanghai University of Medicine $\&$ Health Sciences,   Shanghai 200237, PR China}

\begin{abstract}The aim of this paper is to give new characterizations of some fundamental issues about idempotents. In the general setting of adjointable operators on Hilbert $C^*$-modules,
a new term of quasi-projection pair is introduced. For each   idempotent $Q$, a projection $m(Q)$, called the matched projection of $Q$, is constructed. It is shown that
$Q$ and $m(Q)$ as idempotents are homotopic, and  $\big(m(Q),Q\big)$ is a quasi-projection pair.
Some formulas for $m(Q)$ are derived. Based on these formulas,  representations and norm estimations associated with $m(Q)$ are dealt with.
\end{abstract}


\end{frontmatter}



\section{Introduction}\label{sec:Intro}

A huge number of papers focusing on projections or  idempotents
can be found in the literature (see e.\,g.\,\cite{Afriat,Ando02,Aurel,Blackadar-02,Bottcher-Spitkovsky,Corach,Halmos,Halpern,Kaftal,Mesland,Paulsen,Pedersen02,Raeburn,Roch}). One of the main topics is the study of the relationship between projections and idempotents. Given an idempotent $Q$ on a Hilbert space $H$ (or more generally on a Hilbert $C^*$-module $H$ \cite{Lance,MT,Paschke}), let
$P_{\mathcal{R}(Q)}$ and $P_{\mathcal{N}(Q)}$ be the projections from $H$ onto
the range and the null space of $Q$, respectively. Some well-known results clarifying the relationship among $Q$, $P_{\mathcal{R}(Q)}$ and $P_{\mathcal{N}(Q)}$ can be found in the literature. For instance, it is shown in \cite[Theorem~4.2]{Afriat}  that $$Q=(I-P_{\mathcal{R}(Q)}P_{\mathcal{N}(Q)})^{-1}P_{\mathcal{R}(Q)} (I-P_{\mathcal{R}(Q)}P_{\mathcal{N}(Q)}),$$ where
$I$ denotes the identity operator on $H$.
Another formula for $Q$ in terms of $P_{\mathcal{R}(Q)}$ and $P_{\mathcal{N}(Q)}$ can be found in \cite[Theorem~3.8]{Ando02}, which indicates that
$$Q=P_{\mathcal{R}(Q)}\big(P_{\mathcal{R}(Q)}+\lambda P_{\mathcal{N}(Q)}\big)^{-1}$$
for every non-zero complex number $\lambda$. On the other hand, the projections $P_{\mathcal{R}(Q)}$ and $P_{\mathcal{N}(Q)}$ can be formulated in terms of $Q$ by
\begin{align}
   \label{range projection Q}& P_{\mathcal{R}(Q)}=Q(Q+Q^*-I)^{-1},  \\
   \label{equ:relation of Q and P NQ}& P_{\mathcal{N}(Q)}=(Q-I)(Q+Q^*-I)^{-1}.
\end{align}
The reader is referred to \cite[Theorem~1.3]{Koliha} for a detailed proof of \eqref{range projection Q}, which gives
 \eqref{equ:relation of Q and P NQ} by a simple replacement of $Q$ with $I-Q$.  There are also many norm equalities and norm inequalities
related to $Q$, $P_{\mathcal{R}(Q)}$ and $P_{\mathcal{N}(Q)}$ (see e.\,g.\,\cite[Theorem~3.7]{Ando02}, \cite[Theorem~7.91]{Aurel}, \cite[Section~5]{Bottcher-Spitkovsky}, \cite[Propositions~3.3 and 4.3]{Koliha}).

Let $Q$ be an idempotent acting on a Hilbert $C^*$-module $H$. Since $\mathcal{R}(Q)$ is closed in $H$, we see that $H$ can be decomposed orthogonally by
 \begin{align*}H=\mathcal{R}(Q)\dotplus \mathcal{N}(Q^*).
 \end{align*}
 So, as in the Hilbert space case, up to unitary equivalence $P_{\mathcal{R}(Q)}$ and $Q$ can be represented by $2\times 2$ block matrices
\begin{equation*}
  P_{\mathcal{R}(Q)} = \left( \begin{array}{cc}  I & 0 \\  0 & 0 \\ \end{array} \right),\quad
  Q = \left(  \begin{array}{cc}I & A \\   0 & 0 \\  \end{array}  \right).
\end{equation*}
 Hence, $P_{\mathcal{R}(Q)}$ and $Q$ as idempotents are homotopic via the  norm-continuous path
 $\begin{pmatrix}                                               I & tA \\
                                               0 & 0 \\
                                             \end{pmatrix}$ for $t\in [0,1]$. Meanwhile, from \eqref{range projection Q} it is easily seen  that the mapping from
 $Q$ to $P_{\mathcal{R}(Q)}$ is norm-continuous. Due to the above homotopic equivalence and  norm continuity,
 the projection $P_{\mathcal{R}(Q)}$ plays a key role in the study of $K$-theory for $C^*$-algebras \cite{Blackadar,RLL,Wegge-Olsen}.

Although much progress has been made in the study of a general  idempotent $Q$ on a Hilbert $C^*$-module $H$, there are still some fundamental issues that are not yet clear. Part of which can be listed as follows:

\begin{problem}\label{prob-01} Is it possible to derive new formulas for $Q$ based on a projection $P$ on $H$ such that $P\ne P_{\mathcal{R}(Q)}$ and $P$  has the same two properties of $P_{\mathcal{R}(Q)}$ mentioned as above?
\end{problem}

\begin{problem}\label{prob-02} Is it true that $\|P_{\mathcal{R}(Q)}-Q\|\le \|P-Q\|$ for every projection $P$ on $H$?
\end{problem}

 In this paper we will introduce a new term of quasi-projection pair $(P,Q)$, in which $P$ is a projection, while $Q$ is an idempotent satisfying certain conditions stated as \eqref{conditions 1 for qpp} and \eqref{conditions 2 for qpp}. For every non-projection idempotent $Q$ on a Hilbert $C^*$-module $H$,  we will show in Theorem~\ref{existence thm of quasi-projection pair} that
there always exist a projection denoted  by $m(Q)$ and an invertible operator $W$ on $H$ such that $\big(m(Q),Q\big)$ is a quasi-projection pair satisfying \eqref{new decomposition for an arbitrary idempotent}. Thus, a new formula for $Q$ is derived  based on $m(Q)$, which is called the matched projection of $Q$. It follows easily from \eqref{new decomposition for an arbitrary idempotent} that
 $m(Q)$ and $Q$ as idempotents are homotopy equivalent. Furthermore,  it follows directly  from Theorem~\ref{thm:formula for the matched projection} that the mapping from $Q$ to $m(Q)$ is norm-continuous. Hence, a positive answer to Problem~\ref{prob-01} is given.

 To answer the second problem, we first investigate the special case that $Q$ is a non-projection idempotent in $M_2(\mathbb{C})$. In this case, we have managed to show that $\|m(Q)-Q\|\le \|P-Q\|$ for every projection $P\in M_2(\mathbb{C})$. Then we move to the general case of Hilbert $C^*$-module operators. We will show in Theorem~\ref{thm:m Q is partially best} that the same inequality is also true whenever $(P,Q)$ is a quasi-projection pair. In addition, we will show that $\|m(Q)-Q\|<\|P_{\mathcal{R}(Q)}-Q\|$ for every non-projection idempotent $Q$ on a Hilbert $C^*$-module (see Theorem~\ref{thm:significance} for the details). Thus, a negative answer to Problem~\ref{prob-02} is provided.

The matched projection $m(Q)$ introduced in this paper for a general idempotent $Q$ on a Hilbert $C^*$-module is a brand new object, which has many interesting features. Based on the formula for $m(Q)$ newly derived in Theorem~\ref{thm:formula for the matched projection},
some fundamental results on $m(Q)$ are obtained in Section~\ref{sec:matched projections}.
For instance, it is proved that $m(Q^*)=m(Q)$  and $m(I-Q)=I-m(Q)$  for every idempotent $Q$ on a Hilbert $C^*$-module (see Theorems~\ref{thm:same matched projection for Q and Q-sta} and \ref{thm:matched projection for I-Q}, respectively).
Some issues on $m(Q)$ (such as that mentioned in Problems~\ref{problem m Q1 Q2} and \ref{problem m Q Q P}) seem very challenging and worth further investigation.

We recall some basic knowledge about adjointable operators on  Hilbert $C^*$-modules. For more details, the reader is referred to
\cite{Lance,MT,Paschke}. Hilbert $C^*$-modules are generalizations of Hilbert spaces by allowing the inner product to take values in a $C^*$-algebra rather than in the complex field. Given Hilbert $C^*$-modules $H$ and $K$ over a $C^*$-algebra $\mathfrak{A}$, let $\mathcal{L}(H,K)$ denote the set of all adjointable operators from $H$ to $K$.
It is known that every adjointable operator $A\in\mathcal{L}(H,K)$ is a bounded linear operator, which is also $\mathfrak{A}$-linear in the sense that
\begin{equation*}\label{equ:keep module operator}A(xa)=A(x)a,\quad\forall\,x\in H, a\in\mathfrak{A}.
\end{equation*}
Let $\mathcal{R}(A)$ and $\mathcal{N}(A)$ denote the range and the null space of $A$, respectively.
  When $H=K$, the notation $\mathcal{L}(H,H)$ is simplified to  $\mathcal{L}(H)$.  Let $\mathcal{L}(H)_+$ denote the set of all positive elements in the $C^*$-algebra $\mathcal{L}(H)$. Given $A\in\mathcal{L}(H)$ and
    $M\subseteq H$, let $\sigma(A)$ denote the spectrum of $A$ with respect to $\mathcal{L}(H)$, and let $A|_M$ denote the restriction of $A$ on $M$. When $A\in\mathcal{L}(H)_+$, the notation $A\ge 0$ is also used to indicate that $A$ is a positive operator on $H$. In the special case that $H$ is a Hilbert space, we use the notation $\mathbb{B}(H)$ instead of $\mathcal{L}(H)$.

Given a unital $C^*$-algebra  $\mathfrak{B}$, its unit is denoted by $I_{\mathfrak{B}}$. When $K$ is a Hilbert module over a $C^*$-algebra $\mathfrak{A}$, the unit of $\mathcal{L}(K)$ will be denoted simply by $I_K$ rather than $I_{\mathcal{L}(K)}$. In most cases of this paper,
 the Hilbert $C^*$-module under consideration is represented by $H$, so the symbol $I_ H$ is furthermore simplified to be $I$. In other words, the italicized letter $I$ is reserved for the identity operator on the Hilbert $C^*$-module  $H$.

Suppose that $H$ is a Hilbert module over a $C^*$-algebra $\mathfrak{A}$. An operator $P\in \mathcal{L}(H)$ is called an idempotent if
$P^2=P$. If furthermore $P^*=P$, then $P$ is called a projection, which is said to be non-trivial if  $P\ne I$ and $P\ne 0$.
A closed submodule $M$ of  $H$ is said to be
orthogonally complemented  in $H$ if $H=M+ M^\bot$, where
$M^\bot=\big\{x\in H:\langle x,y\rangle=0\ \mbox{for every}\ y\in
M\big\}$.
In this case, the projection from $H$ onto $M$ is denoted by $P_M$.
We use the notations  "$\oplus$" and "$\dotplus$" with different meanings. For Hilbert
$\mathfrak{A}$-modules $H_1$ and $H_2$, let $H_1\oplus
H_2=\big\{\binom{h_1}{h_2}:h_i\in H_i, i=1,2\big\}$. If both $H_1$ and $H_2$
are submodules of a Hilbert $\mathfrak{A}$-module such that
$H_1\cap H_2=\{0\}$, then we write $H_1\dotplus
H_2=\{h_1+h_2:h_i\in H_i, i=1,2\}.$

Unless otherwise specified, throughout the rest of this paper $\mathbb{C}$ is the complex field,
$\mathfrak{A}$ is a $C^*$-algebra, $E$, $H$ and $K$ are Hilbert $\mathfrak{A}$-modules, $M_n(\mathfrak{B})$ stands for the $n\times n$ matrix algebra over a $C^*$-algebra $\mathfrak{B}$.

The remaining of this paper is organized as follows. In Section~\ref{sec:quasi-projection pairs}, some basic knowledge about quasi-projection pairs are dealt with. A fundamental problem (Problem~\ref{motivation problem}) is presented at the beginning of this section. Motivated by the elementary solution to this problem, the new term of quasi-projection pair is introduced (see Definition~\ref{defn:quasi-projection pair} for the details). As is shown in Theorem~\ref{thm:four equivalences}, this new object has a nice property of symmetry. Some equivalent conditions for quasi-projection pairs are provided in Theorem~\ref{thm:short description of qpp}, which will play a crucial role in the derivations of Theorem~\ref{thm:same matched projection for Q and Q-sta}, Corollary~\ref{cor:commutativity of P and mq} and Theorem~\ref{thm:m Q is partially best}.
Section~\ref{sec:matched projections} is devoted to the study of the representations for the matched projections.
Inspired by Problem~\ref{motivation problem} and formula \eqref{equ:problem P0} for its solution, the matched projection $m(Q)$ of an idempotent $Q$ is introduced
in Theorem~\ref{existence thm of quasi-projection pair} via its $2\times 2$ block matrix form.
Using the knowledge about the Moore-Penrose inverses of adjointable operators on Hilbert $C^*$-modules, we have managed to figure out a formula for $m(Q)$ in Theorem~\ref{thm:formula for the matched projection}; see \eqref{equ:final exp of P} for the details.  Then
we focus on the applications of this formula for $m(Q)$ to derive some fundamental properties of the matched projections (see Theorems~\ref{thm:invance of matched projection}, \ref{thm:same matched projection for Q and Q-sta} and \ref{thm:matched projection for I-Q}).
This formula for $m(Q)$  is also employed to deal with the ranges of the matched projections (see Remark~\ref{rem:range and kernal of mq}, Theorem~\ref{thm:range of m Q-01}, Corollaries~\ref{cor:range of mq and mqm} and \ref{cor:matched pair H2 and H3 zero}).
Section~\ref{sec:norm equ and inequ} deals mainly with the norm estimations associated with the matched projections.
It is interesting to find out upper bounds of $\|m(Q_1)-m(Q_2)\|$ for two general idempotents $Q_1$ and $Q_2$,
and meanwhile to check the validity of  $\|m(Q_1)-m(Q_2)\|\le \|Q_1-Q_2\|$ for some special idempotents $Q_1$ and $Q_2$. Our main results in these two aspects are
Theorems~\ref{thm:norm relation about mq1-mq2 -01}, \ref{thm:conjecture inequality-special} and \ref{thm:norm up of mq1-mq2 accosiated with two projections}. Given a general idempotent $Q$, it is also interesting to investigate conditions on a projection $P$ such that $\|m(Q)-Q\|\le \|P-Q\|$. As is shown in Theorem~\ref{thm:significance}, such an inequality is always true whenever $P=P_{\mathcal{R}(Q)}$. In fact, some deeper results can be found in Theorems~\ref{thm:significance} and  Remark~\ref{rem:optical of two numbers}. We will prove that this latter inequality is also true when $P=0$, $P=I$, $P=P_{\mathcal{N}(Q)}$ respectively  and in the case that $(P,Q)$ is a quasi-projection pair; see  Corollary~\ref{cor:partial answer of pro.2} and Theorem~\ref{thm:m Q is partially best} for the details. To investigate  norm estimations furthermore, it is significant to derive a concrete formula for $\|m(Q)-Q\|$, which will be fulfilled in Theorem~\ref{thm:norm of mq-q}.
Some applications of Theorem~\ref{thm:norm of mq-q} will be dealt with in the last subsection of this paper; see Remark~\ref{rem: alternative proof} and Theorem~\ref{thm:original NIE} for the details.

\section{Quasi-projection pairs}\label{sec:quasi-projection pairs}

Let $0_2$ and $I_2$ be the zero matrix and the identity matrix in  $M_2(\mathbb{C})$. The subset of $M_2(\mathbb{C})$ consisting of all $2\times 2$ projections is denoted by $\mathcal{P}\big(M_2(\mathbb{C})\big)$.

\begin{problem}\label{motivation problem} Suppose that $Q\in M_2(\mathbb{C})$ is an idempotent which is not a projection, whether there exists $P_0\in \mathcal{P}\big(M_2(\mathbb{C})\big)$ such that $$\|P_0-Q\|=\min\big\{\|P-Q\|:P\in \mathcal{P}\big(M_2(\mathbb{C})\big)\big\}?$$
\end{problem}

To answer the above problem, we notice that up to unitary equivalence,
\begin{equation}\label{idempotent in M two}P_{\mathcal{R}(Q)}=\left(
                       \begin{array}{cc}
                         1 & 0 \\
                         0 & 0 \\
                       \end{array}
                     \right),\quad
Q=\left(
                                         \begin{array}{cc}
                                           1 & a \\
                                           0 & 0 \\
                                         \end{array}
                                       \right)\end{equation}
for some $a\in \mathbb{C}\setminus \{0\}$.
Let $a=|a|e^{i\theta}$, in which $\theta \in [0,2\pi)$. Then
\begin{equation*}\label{equ:2 times 2 complect Q}
Q=U\left(  \begin{array}{cc} 1 & |a| \\ 0 & 0 \\   \end{array}  \right)U^*,
\end{equation*}
 where $U=\mbox{diag}\left(1,e^{-i\theta}\right)$. It is obvious that
$$\big\|P_{\mathcal{R}(Q)}-Q\big\|=|a|< \sqrt{1+|a|^2}=\|I_2-Q\|=\|0_2-Q\|,$$
so we need only to consider non-trivial projections in $M_2(\mathbb{C})$.

Now, let $a\in \mathbb{C}\setminus \{0\}$ be fixed and $P\in \mathcal{P}\big(M_2(\mathbb{C})\big)$ be arbitrary that is non-trivial. Put $P^\prime= U^*PU$. Clearly, $P^\prime\in \mathcal{P}\big(M_2(\mathbb{C})\big)$ which is also non-trivial. So by Halmos' two projections theorem for $P_{\mathcal{R}(Q)}$ and $P^\prime$ \cite[Theorem~2.10]{Xu-Yan}, there exist $z\in \mathbb{C}$ with $|z|=1$ and $t\in [0,\pi]$ such that
\begin{equation*}\label{matrix P prime on M2C}
P^\prime = \left(
   \begin{array}{cc}
     \cos^2 t & \bar{z} |\cos t\sin t|\\
     z|\cos t\sin t|  & \sin^2 t \\
   \end{array}
  \right).
\end{equation*}
Hence,
\begin{equation}\label{matrix P on M2C}
P=UP^\prime U^*=U\left(
   \begin{array}{cc}
     \cos^2 t & \bar{z} |\cos t\sin t|\\
     z|\cos t\sin t|  & \sin^2 t \\
   \end{array}
  \right) U^*,
\end{equation}
which gives  $\|P-Q\|=\|R\|,$
where $$R=\left(
            \begin{array}{cc}
              -\sin^2 t & \bar{z} |\cos t\sin t|-|a|  \\
               z|\cos t\sin t| &  \sin^2 t \\
            \end{array}
          \right). $$
A simple calculation yields $$RR^*=\left(
                 \begin{array}{cc}
                   \sin^2t+|a|^2-2|a| \cdot |\cos t \sin t|x & -|a|\sin^2 t \\
                   -|a|\sin^2 t & \sin^2t \\
                 \end{array}
               \right),$$
where
\begin{equation*}\label{real part x}x=\frac{1}{2}(z+\bar{z})\in [-1,1].\end{equation*}
Let $\lambda_i(x,t)$ $(i=1,2)$ denote the
eigenvalues of $RR^*$.  Since $RR^*$ is positive semi-definite, we have
\begin{align*}\|R\|^2=&\|RR^*\|=\max\{\lambda_1(x,t), \lambda_2(x,t)\}\\
=&\frac{1}{2}\left[2\sin^2t+|a|\left(\mu(x,t)+\sqrt{\mu^2(x,t)+4\sin^4t}\right)\right]
 \end{align*}
for $ t\in [0,\pi]$ and $x\in [-1,1]$, in which $$\mu(x,t)=|a|-2|\cos t\sin t|x.$$
If $t\in [0,\pi]\setminus \{0,\frac{\pi}{2}, \pi\}$, then $|\cos t \sin t|>0$, so for any such fixed $t$,
$\mu(x,t)$ is monotonically decreasing  with respect to $x\in [-1,1]$, while  the function $u\to u+\sqrt{u^2+4\sin^4t}$ is monotonically increasing on the real line, thus $\mu(x,t) \ge \mu(1,t)$.
The same is true if $t\in \{0,\frac{\pi}{2}, \pi\}$, since in this case $\mu(x,t)\equiv |a|$ for every $x\in [-1,1]$. So it can be concluded that
\begin{align*}
\min_{t\in [0,\pi],\,x\in [-1,1]} \|R\|^2=\min_{t\in [0,\pi],x=1} \|R\|^2=\min_{t\in [0,\pi]}f(t),
\end{align*}
where
\begin{align*}
f(t)=&\frac{1}{2}\left[2\sin^2t+|a|\left(|a|-|\sin(2t)|+\sqrt{\big(|a|-|\sin(2t)|\big)^2+4\sin^4t}\right)\right]\\
=&\frac12 \left[ |a|^2+1-g(t)+|a| \sqrt{ |a|^2+ 2-2g(t)}\right],
\end{align*}
in which
$$g(t)=\cos(2t)+|a|\cdot |\sin(2t)|.$$ It is clear that $g(t)\le b$ for all $t\in [0,\pi]$,
where
\begin{equation*}\label{equ:defn of b}
b=\sqrt{1+|a|^2}.
\end{equation*}
Let $\theta_0\in (0,\frac{\pi}{2})$ be chosen such that
\begin{equation}\label{cs theta 0}\sin(\theta_0)=\frac{|a|}{b},\quad \cos(\theta_0)=\frac{1}{b}.\end{equation}
Put $t_0=\frac{\theta_0}{2}$. Then $t_0\in (0,\frac{\pi}{4})$, so $|\sin(2t_0)|=\sin(2t_0)$ and thus $$g(t_0)=b\big[\cos\theta_0\cos(2t_0)+\sin\theta_0\sin(2t_0)\big]=b.$$
It follows that the function $f(t)$  takes the minimum value at the point $t=t_0$.
Therefore, $$\|P_0-Q\|=\min\big\{\|P-Q\|:P\in \mathcal{P}\big(M_2(\mathbb{C})\big)\big\},$$
where $P_0$ is formulated by \eqref{matrix P on M2C} with $z=1$ and $t=t_0$ therein.
Hence, a positive answer to Problem~\ref{motivation problem} is given. Moreover,
from \eqref{cs theta 0} we obtain
\begin{equation*}\label{equ:SC equ}
\sin(t_0)\cos(t_0)=\frac{\sqrt{b^2-1}}{2b}=\frac{|a|}{2b},
\end{equation*}
which gives
\begin{align}\label{equ:problem P0}
P_0=&U \left(
   \begin{array}{cc}
     \cos^2 t_0 &  |\cos t_0\sin t_0|\\
     |\cos t_0\sin t_0|  & \sin^2 t_0 \\
   \end{array}
  \right)U^*
  =\frac{1}{2b}\left(\begin{array}{cc}
        b+1 & a\\
        \bar{a} & b-1 \\
      \end{array}\right).
\end{align}

Based on the derived formulas for $P_0$ and $Q$, it can be verified directly that $P_0$ and $Q$ satisfy the conditions stated in \eqref{conditions 1 for qpp} and \eqref{conditions 2 for qpp} with $P$ therein be replaced by $P_0$.
Inspired by such a pair of $P_0$ and $Q$, we now introduce a new term of quasi-projection pair as follows.

\begin{definition}\label{defn:quasi-projection pair}An ordered pair $(P,Q)$  is called a quasi-projection pair on $H$ if $P\in\mathcal{L}(H)$ is a projection, while $Q\in\mathcal{L}(H)$ is an idempotent such that
\begin{align}\label{conditions 1 for qpp}&PQ^*P=PQP, \quad PQ^*(I-P)=-PQ(I-P),\\
\label{conditions 2 for qpp}&(I-P)Q^*(I-P)=(I-P)Q(I-P).
\end{align}
\end{definition}

Note that if $(P,Q)$ is a quasi-projection pair, then from the second equation in \eqref{conditions 1 for qpp}
we can get $(I-P)Q^*P=-(I-P)QP$ by taking $*$-operation. Also, it is notable that $(P,Q)$ may fail to be a quasi-projection pair even if $P$ and $Q$ are two projections.

\begin{remark}\label{rem:trivial quasi-projection pair} From Definition~\ref{defn:quasi-projection pair}, it is clear that $(P,P)$ is a quasi-projection pair
whenever $P\in\mathcal{L}(H)$ is a projection. In the special case that $P$ is a trivial projection (that is, either $P=I$ or $P=0$), then
for every $Q\in\mathcal{L}(H)$, $(P,Q)$ is a quasi-projection pair if and only if $Q$ is a projection.
\end{remark}

We provide two useful characterizations of quasi-projection pairs.
\begin{theorem}\label{thm:four equivalences} Suppose that $P\in\mathcal{L}(H)$ is a projection and $Q\in\mathcal{L}(H)$ is an idempotent.  If one element in
  \begin{equation*}\label{Sigma eight elements}\Sigma:=\big\{(A,B): A\in \{P, I-P\}, B\in \{Q,Q^*, I-Q,I-Q^*\}\big\}\end{equation*}
  is a  quasi-projection pair, then all the remaining elements in $\Sigma$ are  quasi-projection pairs.
 \end{theorem}
\begin{proof} Without loss of generality, we may assume that $(P,Q)$ is a quasi-projection pair. That is, \eqref{conditions 1 for qpp} and \eqref{conditions 2 for qpp} are  both satisfied. It is easy to verify that \eqref{conditions 1 for qpp} and \eqref{conditions 2 for qpp} are also valid whenever $(P,Q)$ is replaced with $(P,Q^*)$, $(P,I-Q)$ and $(I-P,Q)$, respectively.
This ensures the validity of \eqref{conditions 1 for qpp} and \eqref{conditions 2 for qpp} for all elements in $\Sigma$.
\end{proof}

\begin{theorem}\label{thm:short description of qpp} Suppose that $P\in \mathcal{L}(H)$ is a projection and $Q\in\mathcal{L}(H)$ is an idempotent. Then the following statements are equivalent:
\begin{itemize}
  \item [\rm{(i)}] $(P,Q)$ is a quasi-projection pair;
  \item [\rm{(ii)}]$Q^*=(2P-I)Q(2P-I)$;
   \item [\rm{(iii)}]$|Q^*|=(2P-I)|Q|(2P-I)$.
 \end{itemize}
\end{theorem}
\begin{proof} (i)$\Longrightarrow$(ii). From \eqref{conditions 1 for qpp}--\eqref{conditions 2 for qpp}, we have
\begin{align*}
  Q^* & =PQ^*P+PQ^*(I-P)+(I-P)Q^*P+(I-P)Q^*(I-P) \\
   & =PQP-PQ(I-P)-(I-P)QP+(I-P)Q(I-P)\\
   &=PQ(2P-I)-(I-P)Q(2P-I)=(2P-I)Q(2P-I).
\end{align*}

(ii)$\Longrightarrow$(i). Clearly, $P(2P-I)=(2P-I)P=P$. So
$$PQ^*P=P(2P-I)Q(2P-I)P=PQP.$$
The other two equations in \eqref{conditions 1 for qpp}--\eqref{conditions 2 for qpp} can be derived similarly.

(ii)$\Longrightarrow$(iii). As $2P-I$ is a self-adjoint unitary, it follows directly from $Q^*=(2P-I)Q(2P-I)$ that
\begin{align}\label{mark for qpp-01}
   & Q=(2P-I)Q^*(2P-I),\quad |Q|=(2P-I)|Q^*|(2P-I).
\end{align}

(iii)$\Longrightarrow$(ii). The closedness of $\mathcal{R}(Q)$ in $H$ ensures that $Q$ as well as $Q^*$, $|Q|$ and  $|Q^*|$ is Moore-Penrose invertible \cite[Theorem~2.2]{Xu-sheng}. From the basic knowledge about the Moore-Penrose inverse, we have
\begin{equation}\label{range projection ito MP inverse}|Q^*|^\dag\cdot |Q^*|=|Q^*|\cdot |Q^*|^\dag=P_{\mathcal{R}(Q)},\end{equation}
where $|Q^*|^\dag$ denotes the  Moore-Penrose inverse of $|Q^*|$. In virtue of  $Q^*\cdot |Q^*|^\dag \cdot Q\geq 0$,  and
$$\big(Q^*\cdot |Q^*|^\dag \cdot Q\big)^2=Q^*|Q^*|^\dag\cdot |Q^*|^2\cdot |Q^*|^\dag Q=Q^*\cdot P_{\mathcal{R}(Q)}\cdot Q=|Q|^2,$$
we see that
\begin{equation}\label{equ:another expresion of midQmid}
Q^*\cdot |Q^*|^\dag \cdot Q=|Q|.
\end{equation}
Therefore,
\begin{equation}\label{produces Q}|Q^*|\cdot |Q|=(|Q^*|\cdot Q^*)|Q^*|^\dag \cdot Q=(|Q^*|\cdot |Q^*|^\dag) \cdot Q=Q.
\end{equation}
Exchanging $Q$ with $Q^*$ yields
\begin{equation}\label{produces Q-star}|Q|\cdot |Q^*|=Q^*.
\end{equation}

Suppose now that $|Q^*|=(2P-I)|Q|(2P-I)$, then clearly the second equation in \eqref{mark for qpp-01} is satisfied, so according to \eqref{produces Q} and \eqref{produces Q-star}  we have
\begin{align*}(2P-I)Q(2P-I)=&(2P-I)|Q^*|(2P-I)\cdot (2P-I)|Q|(2P-I)\\
=&|Q|\cdot |Q^*|=Q^*. \qedhere
\end{align*}
\end{proof}

\section{Representations for the matched projections}\label{sec:matched projections}

Let $Q$ be an arbitrary idempotent on $H$. To investigate the existence of a projection $P$ that makes $(P,Q)$ to be a quasi-projection pair, we need some basic knowledge about the $2\times 2$ matrix representation for an operator
induced by a projection.  Given a non-trivial projection $P\in\mathcal{L}(H)$,  let $U_P: H\to \mathcal{R}(P)\oplus \mathcal{N}(P)$ be the unitary defined  by
\begin{equation}\label{equ:unitary operator induced by P}U_Ph=\binom{Ph}{(I-P)h},\quad \forall\, h\in H.
\end{equation}
It is clear that for every $h_1\in\mathcal{R}(P)$ and $h_2\in\mathcal{N}(P)$,
\begin{equation*}\label{expression of inverse of U P}U_P^{*}\binom{h_1}{h_2}=h_1+h_2.\end{equation*}
For every $T\in \mathcal{L}(H)$,  we have  (see e.\,g.\,\cite[Section~2.2]{Xu-Wei-Gu})
 \begin{equation}\label{equ:block matrix T}
U_PTU_P^{*}=\left(
            \begin{array}{cc}
              PTP|_{\mathcal{R}(P)} & PT(I-P)|_{\mathcal{N}(P)} \\
              (I-P)TP|_{\mathcal{R}(P)} &(I-P)T(I-P)|_{\mathcal{N}(P)}\\
            \end{array}
          \right).
\end{equation}
Specifically,
 \begin{equation*}\label{equ:block matrix P}
U_PPU_P^{*}=\left(
            \begin{array}{cc}
              I_{\mathcal{R}(P)} & 0 \\
              0 &0\\
            \end{array}
          \right),\quad U_P(I-P)U_P^{*}=\left(
            \begin{array}{cc}
              0 & 0 \\
              0 &I_{\mathcal{N}(P)}\\
            \end{array}
          \right).
\end{equation*}
Conversely, given every $X=\left(\begin{array}{ccc} X_{11}&
X_{12}\\X_{21}&X_{22}\end{array}\right)\in \mathcal{L}\big(\mathcal{R}(P)\oplus \mathcal{N}(P)\big)$ and every $h\in H$, we have
\begin{align*}U_P^{*}XU_P(h)=&U_P^{*}X\binom{Ph}{h-Ph}=(X_{11}+X_{21})Ph+(X_{12}+X_{22})(I-P)h.
\end{align*}
It follows that
\begin{equation}\label{eqn:reverse direction}U_P^{*}XU_P=(X_{11}+X_{21})P+(X_{12}+X_{22})(I-P).\end{equation}

 Inspired by Problem~\ref{motivation problem} and formula \eqref{equ:problem P0} for its solution, we are now able to provide an interesting result as follows.
\begin{theorem}\label{existence thm of quasi-projection pair} For every idempotent $Q\in\mathcal{L}(H)$,
there exist a projection $P\in\mathcal{L}(H)$ and an invertible operator $W\in\mathcal{L}(H)$ such that $(P,Q)$ is a quasi-projection pair satisfying
\begin{equation}\label{new decomposition for an arbitrary idempotent}Q=W^{-1}PW\quad\mbox{and}\quad  \|I-W\|<1.
\end{equation}
In particular, $P$ and $Q$ as idempotents are homotopy equivalent.
\end{theorem}
\begin{proof}Denote $I_{\mathcal{R}(Q)}$ and $I_{\mathcal{N}(Q^*)}$  simply by $I_1$ and $I_2$, respectively. Also, let $I_1\oplus I_2$ be simplified as $\widetilde{I}$. By Remark~\ref{rem:trivial quasi-projection pair}, we may  as well assume that $Q$ is not a projection. Let $P_{\mathcal{R}(Q)}$ be the projection from $H$ onto $\mathcal{R}(Q)$, which can be formulated by \eqref{range projection Q}.
 Let $U_{P_{\mathcal{R}(Q)}}$ be the unitary operator induced by $P_{\mathcal{R}(Q)}$ (see \eqref{equ:unitary operator induced by P}), and let
 \begin{equation}\label{equ:1st defn of AB}A=\big(Q-P_{\mathcal{R}(Q)}\big)|_{\mathcal{N}(Q^*)},\quad B=(AA^*+I_1)^\frac12.\end{equation}
  Since
\begin{equation*}\label{equ:simplification for idempotent} P_{\mathcal{R}(Q)}Q\big(I-P_{\mathcal{R}(Q)}\big)=Q\big(I-P_{\mathcal{R}(Q)}\big)=Q-P_{\mathcal{R}(Q)},\end{equation*} from \eqref{equ:block matrix T} we obtain
\begin{equation}\label{defn of widetile Q-1}\widetilde{Q}:= U_{P_{\mathcal{R}(Q)}}\cdot Q\cdot U_{P_{\mathcal{R}(Q)}}^{*}=\left(
                                                       \begin{array}{cc}
                                                         I_1 & A \\
                                                         0 & 0 \\
                                                       \end{array}
                                                     \right).\end{equation}
Let $\widetilde{P},\widetilde{W}\in\mathcal{L}\big(\mathcal{R}(Q)\oplus\mathcal{N}(Q^*)\big)$ be defined by
\begin{align}&\label{defn of widetilde Q wrt A B}
\widetilde{P}=\frac12 \left(\begin{array}{cc}
        (B+I_1)B^{-1} &  B^{-1}A\\
        A^*B^{-1} &  A^*\big[B(B+I_1)\big]^{-1}A \\
      \end{array}\right),\\
&\nonumber \widetilde{W}= \frac{1}{2} \left(\begin{array}{cc}
        B^{-1} & 0\\
        A^*\big[B(B+I_1)\big]^{-1} &  2I_2\\
      \end{array}\right).
\end{align}
Obviously, $\widetilde{P}^*=\widetilde{P}$ and by \eqref{equ:1st defn of AB} we have
\begin{align}\label{equ:needed for abs of A star-more} AA^*=B^2-I_1.
\end{align}
 In virtue of  \eqref{defn of widetilde Q wrt A B}, \eqref{defn of widetile Q-1} and \eqref{equ:needed for abs of A star-more},  by direct computations we have
 \begin{align*}&\widetilde{P}^2=\widetilde{P},\quad \widetilde{P}\left(\widetilde{Q}+\widetilde{Q}^*\right)(\widetilde{I} -\widetilde{P})=0,\\
 &\widetilde{P}\widetilde{Q}\widetilde{P}=\frac14\left(
                                                         \begin{array}{cc}
                                                           (B+I_1)^2B^{-1} & (B+I_1)B^{-1}A\\
                                                            A^*(B+I_1)B^{-1}& A^*B^{-1}A\\
                                                         \end{array}
                                                       \right), \\
&(\widetilde{I} -\widetilde{P})\widetilde{Q}(\widetilde{I} -\widetilde{P})=\frac14 \left(
                                                         \begin{array}{cc}
                                                           -(B-I_1)^2B^{-1} & (B-I_1)B^{-1}A\\
                                                            A^*(B-I_1)B^{-1}& -A^*B^{-1}A\\
                                                         \end{array}
                                                       \right).
\end{align*}
This shows that $(\widetilde{P},\widetilde{Q})$ is a quasi-projection pair.

Evidently, $\widetilde{W}$ is an invertible operator in $\mathcal{L}\big(\mathcal{R}(Q)\oplus\mathcal{N}(Q^*)\big)$ whose inverse is given by
$$\widetilde{W}^{-1}=\left(
                       \begin{array}{cc}
                         2B & 0 \\
                         -A^*(B+I_1)^{-1} & I_2 \\
                       \end{array}
                     \right).
$$
Bearing \eqref{equ:needed for abs of A star-more} in mind, we have
\begin{align*}\widetilde{W}\widetilde{Q}=\frac{1}{2}\left(
                                           \begin{array}{cc}
                                            B^{-1} & B^{-1}A \\
                                            A^*\big[B(B+I_1)\big]^{-1}  & A^*\big[B(B+I_1)\big]^{-1}A \\
                                           \end{array}
                                         \right)=\widetilde{P}\widetilde{W}.
\end{align*}
This shows that $\widetilde{Q}=\widetilde{W}^{-1} \widetilde{P}\widetilde{W}$. Let
$$P=U_{P_{\mathcal{R}(Q)}}^{*}\cdot \widetilde{P}\cdot U_{P_{\mathcal{R}(Q)}}\quad\mbox{and}\quad  W=U_{P_{\mathcal{R}(Q)}}^{*}\cdot \widetilde{W}\cdot U_{P_{\mathcal{R}(Q)}}.$$
Then clearly, $Q=W^{-1}PW$.

Next, we show that $\|W-I\|<1$ and $P$ and $Q$ are homotopy equivalent. A direct computation yields $$(\widetilde{W}-\widetilde{I})^*(\widetilde{W}-\widetilde{I})=\left(
  \begin{array}{cc}
    S & 0 \\
    0 & 0 \\
  \end{array}
\right),$$
where
\begin{align*}S=&\frac{1}{2}(B^2+B)^{-1}(2B^2-I_1)\le \frac{1}{2}(B^2+B)^{-1}\cdot 2B^2=(B+I_1)^{-1}B.
\end{align*}
Therefore,
$$\big\|\widetilde{W}-\widetilde{I}\big\|^2=\|S\|\leq \big\|(B+I_1)^{-1}B\big\|=\frac{\|B\|}{\|B\|+1}<1.$$
Hence, $\|W-I\|=\big\|\widetilde{W}-\widetilde{I}\big\|<1$. It is a routine matter to put
$$W_t=I+t(W-I)\quad\mbox{and}\quad Q(t)=W_t^{-1}PW_t$$
for $t\in [0,1]$. Obviously, $Q(t)$ is a norm-continuous path of idempotents in $\mathcal{L}(H)$ which starts at $P$ and ends at $Q$.
\end{proof}

\begin{definition}\label{defn:matched projection} Let $Q\in\mathcal{L}(H)$ be an idempotent. If $Q$ is a projection, then $Q$ is referred to be the matched projection of $Q$; otherwise, the projection $P$ constructed in the proof of Theorem~\ref{existence thm of quasi-projection pair} will be  called the matched projection of $Q$. In each case, $(P,Q)$ is called a matched pair, and meanwhile $P$ is denoted by $m(Q)$.
\end{definition}

\begin{remark}Suppose that  $Q\in \mathcal{L}(H)$ is non-projection idempotent. As is shown in
Theorem~\ref{existence thm of quasi-projection pair},  $m(Q)$ and $Q$ as idempotents are homotopy equivalent. Hence, the similarity of
$m(Q)$ and $Q$ is ensured \cite[Proposition~4.3.3]{Blackadar}; see \eqref{new decomposition for an arbitrary idempotent} and Lemma~\ref{lem:pre case p-mq le 1} for the different construction of the associated invertible operators.
It is known (see e.\,g.\,\cite[Proposition~4.6.2]{Blackadar}) that $Q$ is also homotopic to $P_{\mathcal{R}(Q)}$.
The distances of $m(Q)$ and $P_{\mathcal{R}(Q)}$ away from $Q$ will be compared in Theorem~\ref{thm:significance}.
\end{remark}

Throughout the rest of this paper, $m(Q)$ denotes the matched projection of an idempotent $Q$.
It follows directly  from the above definition that  a matched pair  is a quasi-projection pair. In Theorem~\ref{thm:m Q is partially best}, we will show that among all quasi-projection pairs $(P,Q)$, the matched pair $\big(m(Q),Q\big)$ is the best in the sense of the distances between $P$ and $Q$.

Given an idempotent $Q\in\mathcal{L}(H)$, it is meaningful to derive concrete formulas for $m(Q)$. This will be dealt with in our next theorem.

\begin{theorem}\label{thm:formula for the matched projection}For every idempotent $Q\in\mathcal{L}(H)$, we have
\begin{equation}\label{equ:final exp of P}m(Q)=\frac12\big(|Q^*|+Q^*\big)|Q^*|^\dag\big(|Q^*|+I\big)^{-1}\big(|Q^*|+Q\big),
\end{equation}
where $|Q^*|^\dag$ denotes the Moore-Penrose inverse of $|Q^*|$, which can be formulated by
\begin{equation}\label{formula for the MP inverse-01}|Q^*|^\dag=\left(P_{\mathcal{R}(Q)}P_{\mathcal{R}(Q^*)}P_{\mathcal{R}(Q)}\right)^\frac12.
\end{equation}
\end{theorem}
\begin{proof}First we assume that $Q$ is a projection. In this case, the right side of \eqref{equ:final exp of P} turns out to be
$2Q(Q+I)^{-1}$, which is evidently equal to $Q$. So, the conclusion holds.

Suppose now that $Q$ is not a projection.  Utilizing \eqref{defn of widetilde Q wrt A B} and \eqref{eqn:reverse direction} yields
\begin{align*}m(Q)=& \frac12 (B+I_1+A^*)B^{-1}P_{\mathcal{R}(Q)}\\
&+\frac12 \left[B^{-1}A+A^*\big[B(B+I_1)\big]^{-1}A\right]\left(I-P_{\mathcal{R}(Q)}\right),
\end{align*}
where $P_{\mathcal{R}(Q)}$, $A$ and $B$ are given by \eqref{range projection Q} and \eqref{equ:1st defn of AB}, respectively.
Let
\begin{equation}\label{equ:defn of widetilde A and B}\widetilde{A}=Q-P_{\mathcal{R}(Q)},\quad  \widetilde{B}=\big(\widetilde{A}\widetilde{A}^*+I\big)^\frac12,\end{equation}
and let $\mathfrak{\widetilde{B}}$ (resp.\,$\mathfrak{B}$) be the $C^*$-subalgebra of $\mathcal{L}(H)$ \big(resp.\,$\mathcal{L}[\mathcal{R}(Q)]$\big) generated by
$I$ and $\widetilde{A}\widetilde{A}^*$ (resp.\,$I_1$ and $AA^*$). Since $A$ is given by \eqref{equ:1st defn of AB} and $A^*=\widetilde{A}^*|_{\mathcal{R}(Q)}$, we see that
$\widetilde{A}\widetilde{A}^*\xi=AA^*\xi$ for every $\xi\in \mathcal{R}(Q)$, which gives
$p(\widetilde{A}\widetilde{A}^*)\xi=p(AA^*)\xi$, where $p$ is an arbitrary polynomial. So for every function $f$ that is continuous  on $\big[0,\|\widetilde{A}\|^2\big]$, we have
\begin{equation}\label{equ:conti func AA star}
f(AA^*)\xi=f(\widetilde{A}\widetilde{A}^*)\xi,\quad \forall\xi\in \mathcal{R}(Q).
\end{equation}
Particularly, if we take $f(t)=\frac{1}{\sqrt{1+t}}$, then due to \eqref{equ:1st defn of AB}, \eqref{equ:defn of widetilde A and B} and $\mathcal{R}(A)\subseteq \mathcal{R}(Q)$, we may use \eqref{equ:conti func AA star} to get $B^{-1}A=\widetilde{B}^{-1}A$. Similarly,
$$B^{-1}P_{\mathcal{R}(Q)}= \widetilde{B}^{-1}P_{\mathcal{R}(Q)},\quad \big[B(B+I_1)\big]^{-1}A=\big[\widetilde{B}(\widetilde{B}+I)\big]^{-1}A.$$
The equations above together with $A\left(I-P_{\mathcal{R}(Q)}\right)=\widetilde{A}$ yield
\begin{equation*}\label{temp rep for P}
m(Q)=\frac{1}{2}Y+\frac{1}{2}Z,
\end{equation*}
 where
\begin{align}
  \label{temp rep for P Y} & Y=\left(\widetilde{B}+I+\widetilde{A}^*\right)\widetilde{B}^{-1}P_{\mathcal{R}(Q)}, \\
  \label{temp rep for P Z}  & Z=\left[I+\big[(\widetilde{B}+I)^{-1}\widetilde{A}\big]^*\right]\widetilde{B}^{-1}\widetilde{A}.
\end{align}
It is clear that
$$QP_{\mathcal{R}(Q)}=P_{\mathcal{R}(Q)},\quad \left(I-P_{\mathcal{R}(Q)}\right)QQ^*=0,$$
which lead to
$$P_{\mathcal{R}(Q)}Q^*=P_{\mathcal{R}(Q)},\quad QQ^*\left(I-P_{\mathcal{R}(Q)}\right)=0$$
by taking $*$-operation.
Hence,
\begin{align*}\widetilde{B}^2=\widetilde{A}\widetilde{A}^*+I=QQ^*+I-P_{\mathcal{R}(Q)}=\left[|Q^*|+I-P_{\mathcal{R}(Q)}\right]^2.
\end{align*}
Therefore,
\begin{equation}\label{equ:new form for wide Q-01}\widetilde{B}=|Q^*|+I-P_{\mathcal{R}(Q)}.
\end{equation}
So, from \eqref{range projection ito MP inverse} and \eqref{equ:new form for wide Q-01} we can obtain
\begin{align}\label{equ:defn of widetilde inv B}&\widetilde{B}^{-1}=|Q^*|^\dag+I-P_{\mathcal{R}(Q)},\\
\label{equ:defn of widetilde inv B+1}&(\widetilde{B}+I)^{-1}=\big(|Q^*|+I\big)^{-1}P_{\mathcal{R}(Q)}+\frac12 (I-P_{\mathcal{R}(Q)}).
\end{align}
It follows from \eqref{equ:defn of widetilde inv B}, \eqref{range projection ito MP inverse}, \eqref{equ:defn of widetilde A and B} and \eqref{equ:defn of widetilde inv B+1} that
\begin{align*}&\widetilde{B}^{-1} Q=|Q^*|^\dag Q,\quad \widetilde{B}^{-1}P_{\mathcal{R}(Q)}=|Q^*|^\dag P_{\mathcal{R}(Q)}=|Q^*|^\dag,\nonumber\\
&\widetilde{B}^{-1}\widetilde{A}=\widetilde{B}^{-1} Q-\widetilde{B}^{-1}P_{\mathcal{R}(Q)}=|Q^*|^\dag(Q-I),\\
&\big(\widetilde{B}+I\big)^{-1}\widetilde{A}=\big(|Q^*|+I\big)^{-1}\widetilde{A}.\nonumber
\end{align*}
Substituting \eqref{equ:defn of widetilde A and B} and the above equations into \eqref{temp rep for P Y} and \eqref{temp rep for P Z} yields
\begin{align*}Y=&P_{\mathcal{R}(Q)}+\left(Q^*+I-P_{\mathcal{R}(Q)}\right)|Q^*|^\dag=P_{\mathcal{R}(Q)}+Q^*|Q^*|^\dag,\\
Z=&\left[I+\left(Q^*-P_{\mathcal{R}(Q)}\right)\big(|Q^*|+I\big)^{-1}\right]|Q^*|^\dag(Q-I)\\
=&\left(|Q^*|+I+Q^*-P_{\mathcal{R}(Q)}\right)\big(|Q^*|+I\big)^{-1}|Q^*|^\dag(Q-I)\\
=&\left(|Q^*|+I+Q^*-P_{\mathcal{R}(Q)}\right)|Q^*|^\dag\big(|Q^*|+I\big)^{-1}(Q-I)\\
=&\big(|Q^*||Q^*|^\dag+Q^*|Q^*|^\dag\big)\big(|Q^*|+I\big)^{-1}(Q-I)\\
=&\big(P_{\mathcal{R}(Q)}+Q^*|Q^*|^\dag\big)\big(|Q^*|+I\big)^{-1}(Q-I).
\end{align*}
Based on the above expressions for $Y$ and $Z$, we have
\begin{align*}m(Q)=&\frac12\big(P_{\mathcal{R}(Q)}+Q^*|Q^*|^\dag\big)\left[I+\big(|Q^*|+I\big)^{-1}(Q-I)\right]\nonumber\\
=&\frac12\big(P_{\mathcal{R}(Q)}+Q^*|Q^*|^\dag\big)\big(|Q^*|+I\big)^{-1}\left[\big(|Q^*|+I\big)+Q-I\right]\\
=&\frac12\big(P_{\mathcal{R}(Q)}+Q^*|Q^*|^\dag\big)\big(|Q^*|+I\big)^{-1}\big(|Q^*|+Q\big)\\
=&\frac12\big(|Q^*|+Q^*\big)|Q^*|^\dag\big(|Q^*|+I\big)^{-1}\big(|Q^*|+Q\big).
\end{align*}
Moreover, in view of
\begin{align*}
   & Q^\dag=Q^\dag Q\cdot QQ^\dag=P_{\mathcal{R}(Q^*)}P_{\mathcal{R}(Q)}, \\
   & (QQ^*)^\dag=(Q^*)^\dag Q^\dag=(Q^\dag)^* Q^\dag=P_{\mathcal{R}(Q)}P_{\mathcal{R}(Q^*)}P_{\mathcal{R}(Q)},
\end{align*}
we obtain
\begin{equation*}\label{equ:decomposition of Qdag}
|Q^*|^\dag=\big[(QQ^*)^{\frac{1}{2}}\big]^\dag=\big[(QQ^*)^\dag\big]^{\frac{1}{2}}
=\big(P_{\mathcal{R}(Q)}P_{\mathcal{R}(Q^*)}P_{\mathcal{R}(Q)}\big)^{\frac{1}{2}}.
\end{equation*}
Hence, the desired conclusion follows.
\end{proof}

\begin{remark}Combining \eqref{equ:final exp of P} and \eqref{formula for the MP inverse-01} with \eqref{range projection Q}, we see that the map from an idempotent $Q$ to its matched projection $m(Q)$ is norm-continuous (see also Lemma~\ref{lem:pre case p-mq le 1} for the derivation of such a conclusion). So, the same conclusion stated in \cite[Proposition~4.6.3]{Blackadar} can be derived immediately from Theorems~\ref{existence thm of quasi-projection pair} and \ref{thm:formula for the matched projection}.
\end{remark}

Next, we provide several direct applications of Theorem~\ref{thm:formula for the matched projection}.

\begin{theorem}\label{thm:invance of matched projection} Suppose that $X$ is a Hilbert space and $\pi:\mathcal{L}(H)\to\mathbb{B}(X)$ is a unital $C^*$-morphism. Then $m\big(\pi(Q)\big)=\pi\big(m(Q)\big)$ for every idempotent $Q\in \mathcal{L}(H)$.
\end{theorem}
\begin{proof} Since $\pi$ is a unital $C^*$-morphism, by  \eqref{range projection Q} we have
\begin{equation}\label{equ:map range of range projection}\pi\big(P_{\mathcal{R}(Q)}\big)=\pi(Q)\big[\pi(Q)+\pi(Q)^*-I]^{-1}
=P_{\mathcal{R}[\pi(Q)]}.\end{equation}
Replacing $Q$ with $Q^*$ yields
\begin{equation}\label{equ:map range of range projection+2}\pi\big(P_{\mathcal{R}(Q^*)}\big)=P_{\mathcal{R}[\pi(Q^*)]}=P_{\mathcal{R}[\pi(Q)^*]}.\end{equation}
Meanwhile, it is evident that
$$\pi(Q)\pi(Q)^*=\pi(QQ^*)=\pi(|Q^*|^2)=\left[\pi(|Q^*|)\right]^2,$$
which means that $\left|\pi(Q)^*\right|=\pi(|Q^*|)$. Furthermore, the equality
$\left|\pi(Q)^*\right|^\dag=\pi\left[|Q^*|^\dag\right]$ can be derived directly from \eqref{formula for the MP inverse-01},
 \eqref{equ:map range of range projection} and \eqref{equ:map range of range projection+2}. So according to \eqref{equ:final exp of P}, we have
\begin{equation*}\pi\big(m(Q)\big)=\frac12\big(|\pi(Q)^*|
+\pi(Q)^*\big) \left|\pi(Q)^*\right|^\dag\big(|\pi(Q)^*|+I\big)^{-1}\big(|\pi(Q)^*|+\pi(Q)\big),
\end{equation*}
which is exactly the formula for $m\big(\pi(Q)\big)$.
\end{proof}

\begin{theorem}\label{thm:same matched projection for Q and Q-sta} For every idempotent $Q\in\mathcal{L}(H)$, we have
$m(Q^*)=m(Q)$.
\end{theorem}
\begin{proof}As $\big(m(Q),Q\big)$ is a quasi-projection pair, by Theorem~\ref{thm:short description of qpp} we have
\begin{equation}\label{qpp for m Q}Q^*=\left[2m(Q)-I\right]Q\left[2m(Q)-I\right].
\end{equation}
Observe that for every unitary $U\in \mathcal{L}(H)$, it is clear that
\begin{align*}|U^*Q^*U|=U^*|Q^*|U,\quad |U^*Q^*U|^\dag =U^*|Q^*|^\dag U,
\end{align*}
so it can be deduced from \eqref{equ:final exp of P} that $m(U^*QU)=U^*m(Q)U$. Specifically, if we put $U=2m(Q)-I$, then
\begin{equation*}m(Q^*)=\left[2m(Q)-I\right]m(Q)\left[2m(Q)-I\right]=m(Q). \qedhere
\end{equation*}

\end{proof}

\begin{corollary}\label{cor:commutativity of P and mq} For every quasi-projection pair $(P,Q)$, we have $ Pm(Q)=m(Q)P$.
\end{corollary}
\begin{proof} By Theorem~\ref{thm:short description of qpp} we have $Q^*=(2P-I)Q(2P-I)$. So, as is shown in the proof of Theorem~\ref{thm:same matched projection for Q and Q-sta}, we can obtain
\begin{align*}
 Pm(Q^*)=&P(2P-I)m(Q)(2P-I)=Pm(Q)(2P-I)\\
 =&2Pm(Q)P-Pm(Q),
\end{align*}
which leads by $m(Q)=m(Q^*)$ to $Pm(Q)=Pm(Q)P$. Taking $*$-operation yields $m(Q)P=Pm(Q)P$. Hence, $Pm(Q)=m(Q)P$.
\end{proof}

\begin{theorem}\label{thm:2 additional formulas for m q} For every idempotent $Q\in\mathcal{L}(H)$, we have
\begin{equation*}\label{equ:2 additional formulas for m q}m(Q)=TT^\dag=VV^*,\quad P_{\mathcal{R}(Q)}=T^\dag T=V^*V,
\end{equation*}
where
\begin{equation*}T=|Q^*|+Q^*,\quad V=\frac{\sqrt{2}}{2} T(|Q^*|^\dag)^{\frac{1}{2}}(I+|Q^*|)^{-\frac{1}{2}}.
\end{equation*}
\end{theorem}
\begin{proof}It follows directly  from \eqref{equ:final exp of P} that
$TW=VV^*=m(Q)$, where
$$W=\frac{1}{2}|Q^*|^\dag(I+|Q^*|)^{-1}(|Q^*|+Q).$$
As $Q$ is an idempotent and
$$\mathcal{R}\big((|Q^*|^\dag)^{\frac{1}{2}}\big)= \mathcal{R}\big(|Q^*|^\dag\big)=\mathcal{R}(|Q^*|)=\mathcal{R}(Q)=\mathcal{R}\big(P_{\mathcal{R}(Q)}\big),$$
 we have $Q|Q^*|=|Q^*|$, so
$$(|Q^*|+Q)|Q^*|=(I+|Q^*|)|Q^*|=(|Q^*|+Q)Q^*.$$
Therefore,
\begin{equation}\label{equ:pre wt}
T^*T=(|Q^*|+Q)(|Q^*|+Q^*)=2(I+|Q^*|)|Q^*|.
\end{equation}
which gives
\begin{equation*}WT=|Q^*|^\dag \cdot |Q^*|=P_{\mathcal{R}(Q)}.\end{equation*}
Note that $\mathcal{R}(T^*)\subseteq \mathcal{R}(Q)$ and $\mathcal{R}(W)\subseteq \mathcal{R}(Q)$, so we have
$$TWT=TP_{\mathcal{R}(Q)}=(P_{\mathcal{R}(Q)}T^*)^*=T,\quad WTW=P_{\mathcal{R}(Q)}W=W.$$
This shows that $T$ is Moore-Penrose invertible such that $W=T^\dag$. Hence, $m(Q)=TT^\dag$ and $P_{\mathcal{R}(Q)}=T^\dag T$. Moreover, from the definition of $V$ and \eqref{equ:pre wt}, we see that
\begin{align*}
  V^*V & = (I+|Q^*|)^{-\frac{1}{2}}(|Q^*|^\dag)^{\frac{1}{2}}\cdot (I+|Q^*|)|Q^*| \cdot (|Q^*|^\dag)^{\frac{1}{2}} (I+|Q^*|)^{-\frac{1}{2}}\\
   & =|Q^*|^\dag \cdot |Q^*|=P_{\mathcal{R}(Q)}. \qedhere
\end{align*}
\end{proof}

\begin{remark}\label{rem:range and kernal of mq} Following the notations as in the proof of Theorem~\ref{thm:2 additional formulas for m q}, we have
\begin{equation*}
  \mathcal{R}\big[m(Q)\big]=\mathcal{R}(TT^\dag)=\mathcal{R}(T),\quad
   \mathcal{N}\big[m(Q)\big]=\mathcal{N}(T^\dag)=\mathcal{N}(T^*).
\end{equation*}
The equations above together with $m(Q)=m(Q^*)$ yield
\begin{align}\label{equ:range of mq}&\mathcal{R}\big[m(Q)\big]=\mathcal{R}(|Q^*|+Q^*)=\mathcal{R}(|Q|+Q),\\
\label{equ:kernel of mq}&\mathcal{N}\big[m(Q)\big]=\mathcal{N}(|Q^*|+Q)=\mathcal{N}(|Q|+Q^*).
\end{align}
It can be deduced from the proof of Theorem~\ref{thm:2 additional formulas for m q} that
\begin{equation}\label{equ:m qstar}
m(Q)Q^*=m(Q)|Q^*|=\frac{1}{2}(|Q^*|+Q^*).
\end{equation}
Since $m(Q)=m(Q^*)$, replacing $Q$ with $Q^*$ in \eqref{equ:m qstar} we can obtain
\begin{equation}\label{equ:m times q}
m(Q)Q=m(Q)|Q|=\frac{1}{2}(|Q|+Q).
\end{equation}
\end{remark}

\begin{theorem}\label{thm:range of m Q-01} For every idempotent $Q\in\mathcal{L}(H)$, we have
 $$\mathcal{R}\big[m(Q)\big]\subseteq \mathcal{R}(Q+Q^*)=\mathcal{R}(|Q^*|+|Q|),$$ and $\mathcal{R}\big[m(Q)\big]=\mathcal{R}(Q+Q^*)$ if and only if $Q$ is a projection.
\end{theorem}
\begin{proof} From \eqref{equ:m times q}, we know that
\begin{equation}\label{equ:relation of q and 1q1}
\big[2m(Q)-I\big]Q=2\cdot \frac{1}{2}(|Q|+Q)-Q=|Q|.
\end{equation}
Replacing $Q$ with $Q^*$ and utilizing Theorem~\ref{thm:same matched projection for Q and Q-sta}, we obtain
$$\big[2m(Q)-I\big]Q^*=|Q^*|.$$
Hence,
\begin{equation}\label{equ:relation of q qstar and 1q1 1qsatr1}
  |Q^*|+|Q|=\big[2m(Q)-I\big](Q+Q^*).
\end{equation}
Taking $*$-operation yields
$$|Q^*|+|Q|=(Q+Q^*)\big[2m(Q)-I\big],$$ which
gives $\mathcal{R}(Q+Q^*)=\mathcal{R}(|Q^*|+|Q|)$, since $2m(Q)-I$ is a unitary.
Let $$Y=(|Q^*|+Q)\big(|Q^*|^\dag\big)^{\frac{1}{2}}.$$ By \eqref{equ:another expresion of midQmid} and \eqref{equ:relation of q qstar and 1q1 1qsatr1}, we have
\begin{align*}
  YY^* & =(|Q^*|+Q^*)|Q^*|^\dag(|Q^*|+Q) \\
   & =|Q^*|+P_{\mathcal{R}(Q)}Q+Q^*P_{\mathcal{R}(Q)}+Q^*|Q^*|^\dag Q\\
   &= |Q^*|+|Q|+Q+Q^*=2(Q+Q^*)m(Q).
\end{align*}
Let $V$ be defined as in Theorem~\ref{thm:2 additional formulas for m q}. Since
$V=Y\cdot \frac{\sqrt{2}}{2} (I+|Q^*|)^{-\frac{1}{2}}$, it is clear that $\mathcal{R}(V)=\mathcal{R}(Y)$, so $\mathcal{R}(Y)$ is closed in $H$ (as $\mathcal{R}(V)$ is). Thus, $\mathcal{R}(Y)=\mathcal{R}(YY^*)$ \cite[Theorem 3.2 and Proposition 3.7]{Lance}. It follows from  Theorem~\ref{thm:2 additional formulas for m q} that
$$\mathcal{R}\big[m(Q)\big]=\mathcal{R}(V)=\mathcal{R}(YY^*)=\mathcal{R}\big[(Q+Q^*)m(Q)\big]\subseteq \mathcal{R}(Q+Q^*).$$

Suppose that $Q\in\mathcal{L}(H)$ is a projection, then $m(Q)=Q^*=Q$, which obviously yields
$\mathcal{R}\big[m(Q)\big]=\mathcal{R}(Q+Q^*)$. Conversely, suppose that $\mathcal{R}\big[m(Q)\big]=\mathcal{R}(Q+Q^*)$, then
$m(Q)(Q+Q^*)=Q+Q^*$. So, \eqref{equ:relation of q qstar and 1q1 1qsatr1} can be simplified as
\begin{equation*}\label{equ:2th relation of q qstar and 1q1 1qsatr1}
  |Q^*|+|Q|=Q+Q^*.
\end{equation*}
The equation above together with $(I-Q)|Q^*|=0$, yields
\begin{align*}(I-Q)|Q|(I-Q^*)=&(I-Q)(|Q^*|+|Q|)(I-Q^*)\\
=&(I-Q)(Q+Q^*)(I-Q^*)=0.
\end{align*}
Hence, $|Q|(I-Q^*)=0$, or equivalently $Q(I-Q^*)=0$, which obviously gives $Q=Q^*$.
\end{proof}

Next, we provide a technical lemma as follows.
\begin{lemma}\label{thm:mQm geq m} Let $Q\in\mathcal{L}(H)$ be an idempotent. Then $m(Q)Q m(Q)\geq m(Q)$.
\end{lemma}
\begin{proof}
From \eqref{formula for the MP inverse-01} we have $\big\| |Q^*|^\dag\big\|\le 1$, so
$$|Q^*|^\dag= P_{\mathcal{R}(Q)} \cdot |Q^*|^\dag \cdot P_{\mathcal{R}(Q)} \le P_{\mathcal{R}(Q)} \cdot \big\| |Q^*|^\dag\big\| \cdot P_{\mathcal{R}(Q)}\le P_{\mathcal{R}(Q)}.$$
Therefore, $$\big(I+|Q^*|\big)^{\frac{1}{2}}|Q^*|^\dag \big(I+|Q^*|\big)^{\frac{1}{2}}=|Q^*|^\dag(I+|Q^*|)=|Q^*|^\dag +P_{\mathcal{R}(Q)}\ge 2|Q^*|^\dag,$$
which leads to
\begin{equation*}\label{equ:Qdag ge Q Qinv}
|Q^*|^\dag(I+|Q^*|)^{-1}= \big(I+|Q^*|\big)^{-\frac{1}{2}}|Q^*|^\dag \big(I+|Q^*|\big)^{-\frac{1}{2}} \leq \frac{1}{2} |Q^*|^\dag.
\end{equation*}
Combining the above inequality with \eqref{equ:final exp of P}, \eqref{equ:another expresion of midQmid}, \eqref{equ:m times q} and \eqref{equ:m qstar} yields
\begin{align*}
  m(Q) & \leq \frac{1}{4}(|Q^*|+Q^*)|Q^*|^\dag(|Q^*|+Q)=\frac{1}{4}(P_{\mathcal{R}(Q)}+Q^*|Q^*|^\dag)\cdot Q(|Q^*|+Q) \\
   & = \frac{1}{4}\big(Q+Q^*|Q^*|^\dag Q\big)(|Q^*|+Q)=\frac{1}{4}(|Q|+Q)(|Q^*|+Q)\\
   &=m(Q)Q\cdot \big[m(Q)Q^*\big]^*=m(Q) Q m(Q).
\end{align*}
This completes the proof.
\end{proof}

As a supplement to Theorem~\ref{thm:range of m Q-01}, we provide a corollary as follows.
\begin{corollary}\label{cor:range of mq and mqm} For every idempotent $Q\in\mathcal{L}(H)$, we have
\begin{equation*}\label{equ:range of mq and mqm}\mathcal{R}\big[m(Q)\big]=\mathcal{R}\big(|Q^*|+|Q|+Q+Q^*\big)=\mathcal{R}\left[\big(m(Q)Q m(Q)\big)^\frac1n\right]\end{equation*}
for every $n\in\mathbb{N}$.
\end{corollary}
\begin{proof}Let $U_{m(Q)}$ be defined by \eqref{equ:unitary operator induced by P} with $P$ therein be replaced by $m(Q)$.
It can be concluded from \eqref{equ:block matrix T} and Lemma~\ref{thm:mQm geq m}  that
\begin{equation*}U_{m(Q)}\left[m(Q)Q m(Q)\right]U_{m(Q)}^*=\left(
                                                         \begin{array}{cc}
                                                           A & 0\\
                                                           0 & 0 \\
                                                         \end{array}
                                                       \right),\end{equation*}
for some $A\in\mathcal{L}(H_1)$ with $H_1=\mathcal{R}\big(m(Q)\big)$ such that $A\ge I_{H_1}$. Therefore,
\begin{equation}\label{blocked matrix m q q m q}U_{m(Q)}\left[\big(m(Q)Q m(Q)\big)^\frac1n\right]U_{m(Q)}^*=\left(
                                                         \begin{array}{cc}
                                                           A^{\frac1n} & 0\\
                                                           0 & 0 \\
                                                         \end{array}
                                                       \right)\quad (\forall\,n\in\mathbb{N}),\end{equation}
which means that
\begin{equation}\label{equ:the lim of mqm}\lim_{n\to\infty}\left\|\big(m(Q)Q m(Q)\big)^\frac1n-m(Q)\right\|=\lim_{n\to\infty}\|A^{\frac1n}-I_{H_1}\|=0,
\end{equation}
and $\mathcal{R}\left[\big(m(Q)Q m(Q)\big)^\frac1n\right]=\mathcal{R}\big[m(Q)\big]$ for every $n\in\mathbb{N}$.
Specifically,  we have
 $\mathcal{R}\big[m(Q)\big]=\mathcal{R}\left[m(Q)Q m(Q)\right]$. By the proof of Lemma~\ref{thm:mQm geq m}, we know that
 $$m(Q)Q m(Q)=\frac{1}{4}(|Q|+Q)(|Q^*|+Q).$$
Therefore, we may use \eqref{produces Q-star} to obtain
\begin{equation*}m(Q)Q m(Q)=\frac14\big[Q^*+|Q|+|Q^*|+Q\big].
\end{equation*}
The desired conclusion follows.
\end{proof}
\begin{theorem}\label{thm:matched projection for I-Q} For every idempotent $Q\in\mathcal{L}(H)$, we have $m(I-Q)=I-m(Q)$.
\end{theorem}
\begin{proof} We begin with some observations concerned with general  idempotents. It is clear that
\begin{equation*}|I-Q^*|^2\cdot |Q|^2=(I-Q)(I-Q^*)\cdot Q^*Q=0.
\end{equation*}
 So we may use  $\mathcal{N}\left(|I-Q^*|^2\right)=\mathcal{N}\left(|I-Q^*|\right)$ and $\mathcal{R}(|Q|^2)=\mathcal{R}(|Q|)=\mathcal{R}(Q^*)$ to get
\begin{equation}\label{equ:relation of q and 1-qstar-01} |I-Q^*|\cdot |Q|=|I-Q^*|\cdot Q^*=0,\quad Q^*|Q|=|Q|.\end{equation}
Observe that
$$ |I-Q^*|^2\cdot (Q^*-Q)=Q^*Q-QQ^*Q=(Q^*-Q)\cdot |Q|^2,$$
which implies that for any polynomial $p$, it holds that
$$p\left(|I-Q^*|^2\right)\cdot (Q^*-Q)=(Q^*-Q)\cdot p\left(|Q|^2\right).$$
By the functional calculus on the commutative $C^*$-subalgebras of $\mathcal{L}(H)$ generated by
$|I-Q^*|^2$ and $|Q|^2$ respectively, it can be deduced that
$$f\left(|I-Q^*|^2\right)\cdot (Q^*-Q)=(Q^*-Q)\cdot f\left(|Q|^2\right)$$
for every function $f$ that is continuous on $[0,+\infty)$. Particularly, if we take $f(t)=\sqrt{t}$, then
we arrive at
\begin{equation*}\label{equ:relation of q and 1-qstar-02}
|I-Q^*|\cdot (Q^*-Q)=(Q^*-Q)\cdot |Q|.
\end{equation*}
A combination of the above equation with \eqref{equ:relation of q and 1-qstar-01} gives
\begin{equation}\label{equ:relation 2 of q and 1-qstar}
  -|I-Q^*|Q=(I-Q)|Q|.
\end{equation}

Now, we are ready to prove the desired conclusion. In virtue of \eqref{equ:range of mq}, we have
$$\mathcal{N}\big[I-m(Q)\big]=\mathcal{R}\big[m(Q)\big]=\mathcal{R}(|Q|+Q).$$
From \eqref{equ:kernel of mq} with $Q$ therein be replaced by $I-Q$, it can be derived that
$$\mathcal{N}\big[m(I-Q)\big]=\mathcal{N}\big[|I-Q^*|+I-Q\big].$$
As both $m(I-Q)$ and $I-m(Q)$ are projections, it needs only to show that $$\mathcal{N}\big[|I-Q^*|+I-Q\big]=\mathcal{R}(|Q|+Q),$$
which can furthermore be reduced to
\begin{equation}\label{one direction inclusion-01}\mathcal{N}\big[|I-Q^*|+I-Q\big]\subseteq \mathcal{R}(|Q|+Q),\end{equation}
since the reverse inclusion is immediate from \eqref{equ:relation of q and 1-qstar-01} and \eqref{equ:relation 2 of q and 1-qstar}.

For each $x\in \mathcal{N}\big[|I-Q^*|+I-Q\big]$, we have
\begin{equation}\label{equ:elments on ker 1-mq}
(I-Q)x=-|I-Q^*|x.
\end{equation}
Let $u=(I+|Q|)^{-1}x$. Then by \eqref{equ:relation 2 of q and 1-qstar} and \eqref{equ:relation of q and 1-qstar-01}, we have
\begin{align}
   &(I-Q)x=(I-Q)(I+|Q|)u=(I-Q)u-|I-Q^*|Qu,  \nonumber\\
 \label{equ:1-qstar x equ 1-qstar u} & |I-Q^*|x=|I-Q^*|(I+|Q|)u=|I-Q^*|u.
\end{align}
Combining the above equations with \eqref{equ:elments on ker 1-mq} gives
$$(I-Q)u=|I-Q^*|Qu-|I-Q^*|u=-|I-Q^*|\cdot(I-Q)u.$$
Therefore, $$\big[I+|I-Q^*|\big]\cdot(I-Q)u=0,$$
which leads by the invertibility of $I+|I-Q^*|$ to $u=Qu$. So, by  \eqref{equ:1-qstar x equ 1-qstar u} and \eqref{equ:relation 2 of q and 1-qstar} we have
$$|I-Q^*|x=|I-Q^*|u=|I-Q^*|Qu=-(I-Q)|Q|u.$$
It follows from \eqref{equ:elments on ker 1-mq} that
\begin{align*}
  x &=Qx+(I-Q)x=Qx-|I-Q^*|x \\
   & =Q(I+|Q|)u+(I-Q)|Q|u=(|Q|+Q)u.
\end{align*}
This shows the validity of \eqref{one direction inclusion-01}.
\end{proof}

\begin{corollary}\label{cor:matched pair H2 and H3 zero}Let $Q\in \mathcal{L}(H)$ be an idempotent. Then
\begin{equation}\label{two intersections are zero}\mathcal{R}\big[m(Q)\big]\cap \mathcal{N}(Q)=\{0\},\quad \mathcal{N}\big[m(Q)\big]\cap \mathcal{R}(Q)=\{0\}.\end{equation}
\end{corollary}
\begin{proof}  Suppose that $x\in \mathcal{N}\big[m(Q)\big]\cap \mathcal{R}(Q)$. Then
$x=Qx$, and by \eqref{equ:kernel of mq} we have $(|Q^*|+Q)x=0$. It follows that $(|Q^*|+I)x=0$, which clearly gives $x=0$.
This shows the validity of the second equation in \eqref{two intersections are zero}. In virtue of $m(I-Q)=I-m(Q)$,
the first equation in \eqref{two intersections are zero} follows immediately.
\end{proof}

\section{Norm estimations associated with the matched projections}\label{sec:norm equ and inequ}

\subsection{Norm upper bounds of $m(Q_1)-m(Q_2)$}

For every projections $P,Q\in\mathcal{L}(H)$, by \cite[Lemma~4.1]{XY2} we have
\begin{equation}\label{equ:KKM equality}\|P-Q\|=\max\big\{\|P(I-Q)\|, \|(I-P)Q\|\big\},
\end{equation}
which is known as Krein-Krasnoselskii-Milman equality in the Hilbert space case. Let $Q_1,Q_2\in\mathcal{L}(H)$ be two arbitrary idempotents. Based on Krein-Krasnoselskii-Milman equality, it can be shown (see the proof of \cite[Proposition~3.3]{Koliha}, for example) that
$$\|P_{\mathcal{R}(Q_1)}-P_{\mathcal{R}(Q_2)}\|\le \|Q_1-Q_2\|.$$
It seems interesting to give answers to the following problem.
\begin{problem}\label{problem m Q1 Q2} Under what conditions it is true that $\|m(Q_1)-m(Q_2)\|\leq \|Q_1-Q_2\|$ for two idempotents $Q_1$ and $Q_2$ in $\mathcal{L}(H)$?\end{problem}

Throughout the rest of this paper, $m(Q), m(Q_1)$ and $m(Q_2)$ are the matched projections of idempotents  $Q,Q_1$ and $Q_2$, respectively. The main purpose of this subsection is to derive some norm upper bounds of $m(Q_1)-m(Q_2)$.
Our first result in this direction reads as follows.
\begin{theorem}\label{thm:norm relation about mq1-mq2 -01} For every idempotents $Q_1,Q_2\in\mathcal{L}(H)$, we have
 \begin{equation*}\label{equ:relation of mq1 mq2}
  \|m(Q_1)-m(Q_2)\|=\inf_{m,n\in \mathbb{N}}\alpha_{m,n}=\inf_{m\in\mathbb{N}} \beta_m=\inf_{n\in\mathbb{N}} \gamma_n,
\end{equation*}
where
\begin{align*}&\alpha_{m,n}=\left\|\left[m(Q_1)Q_1m(Q_1)\right]^\frac1m-\left[m(Q_2)Q_2m(Q_2)\right]^\frac1n\right\|,\\
&\beta_m=\left\|\left[m(Q_1)Q_1m(Q_1)\right]^\frac1m-m(Q_2)\right\|,\\
&\gamma_n=\left\|m(Q_1)-\left[m(Q_2)Q_2m(Q_2)\right]^\frac1n)\right\|.
\end{align*}
  \end{theorem}
\begin{proof} For each idempotent $Q\in\mathcal{L}(H)$ and $n\in\mathbb{N}$, it can be deduced from \eqref{blocked matrix m q q m q} that
$$m(Q)=\big(m(Q)Qm(Q)\big)^\frac1n\cdot \left[\big(m(Q)Qm(Q)\big)^\frac1n\right]^\dag$$
such that
$$\left\|\left[\big(m(Q)Qm(Q)\big)^\frac1n\right]^\dag\right\|=\big\|A^{-\frac1n}\big\|\le 1,$$
since as is observed in the proof of Corollary~\ref{cor:range of mq and mqm} that $A\ge I_{H_1}$.
So if we  put
 $$V_1=m(Q_1)Q_1m(Q_1),\quad V_2=m(Q_2)Q_2m(Q_2),$$
then for every $m,n\in \mathbb{N}$, we have
\begin{align*}
  &\big[I-m(Q_2)\big]m(Q_1)=\big[I-m(Q_2)\big]\big(V_1^\frac{1}{m}-V_2^\frac{1}{n}\big)(V_1^\frac{1}{m})^\dag.
\end{align*}
 Therefore,
\begin{equation*}\label{equ:pre mq1-mq2 01}
\big\|\big[I-m(Q_2)\big]m(Q_1)\big\|\leq \big\|V_1^\frac{1}{m}-V_2^\frac{1}{n}\big\|=\alpha_{m,n},\quad \forall\, m,n\in \mathbb{N}.
\end{equation*}
Exchanging $Q_1$ with $Q_2$ gives
\begin{equation*}\label{equ:pre mq1-mq2 02}
\big\|\big[I-m(Q_1)\big] m(Q_2)\big\|\leq \alpha_{m,n},\quad \forall \, m,n\in \mathbb{N}.
\end{equation*}
So the Krein-Krasnoselskii-Milman equality (see \eqref{equ:KKM equality}) ensures that
\begin{align*}
  \|m(Q_1)-m(Q_2)\| \leq \inf_{m,n\in \mathbb{N}}\alpha_{m,n}.
  \end{align*}
Furthermore, a direct use of \eqref{equ:the lim of mqm} yields
\begin{align*}&\beta_m=\lim_{n\to\infty}\alpha_{m,n},\quad  \gamma_n=\lim_{m\to\infty}\alpha_{m,n},\\
&\lim_{m\to\infty} \beta_m=\lim_{n\to\infty}\gamma_n=\|m(Q_1)-m(Q_2)\|.
\end{align*}
Hence, the desired conclusion follows.
\end{proof}

\begin{remark} Let $Q_1,Q_2\in \mathcal{L}(H)$ be idempotents. Following the notations as in the proof of Theorem~\ref{thm:norm relation about mq1-mq2 -01}, we have $V_1^\dag=m(Q_1)V_1^\dag$, hence $SV_1^\dag=0$, where
$S=m(Q_1)Q_1\big[I-m(Q_1)\big]$. This observation together with
$$\big[I-m(Q_2)\big] m(Q_2)Q_2=0$$ yields
\begin{align*}
  \big[I-m(Q_2)\big]m(Q_1) &  =\big[I-m(Q_2)\big]\big[V_1+S-m(Q_2)Q_2\big]V_1^\dag\\
  &=\big[I-m(Q_2)\big]\big[m(Q_1)Q_1-m(Q_2)Q_2\big]V_1^\dag.
\end{align*}
Consequently,
$$\big\|\big[I-m(Q_2)\big]m(Q_1)\big\|\leq \big\|m(Q_1)Q_1-m(Q_2)Q_2\big\|.$$
So, it can be concluded that
$$ \|m(Q_1)-m(Q_2)\| \leq \big\|m(Q_1)Q_1-m(Q_2)Q_2\big\|$$
by following the line as in the proof of Theorem~\ref{thm:norm relation about mq1-mq2 -01}.
Moreover, if we replace $Q_i$ with $Q_i^*$ for $i=1,2$, then we may use Theorem~\ref{thm:same matched projection for Q and Q-sta} to obtain
$$\|m(Q_1)-m(Q_2)\| \leq \big\|Q_1m(Q_1)-Q_2m(Q_2)\big\|.$$
\end{remark}

To get a partial answer to Problem~\ref{problem m Q1 Q2}, we need two useful lemmas as follows.

\begin{lemma}\label{lem:norm equation about matched projection} For every idempotent $Q\in \mathcal{L}(H)$, we have
\begin{equation*}\label{equ:norm equation about matched projection}
  \big\|m(Q)(I-Q)m(Q)\big\|=\big\|\big[I-m(Q)\big]Q\big[I-m(Q)\big]\big\|.
\end{equation*}
\end{lemma}
\begin{proof} It needs to prove that $\|T\|=\|S\|$, where
\begin{align}\label{equ:auxiliary positive operators}
   & T=-m(Q)(I-Q)m(Q),\quad S=-\big[I-m(Q)\big]Q\big[I-m(Q)\big].
\end{align}
In virtue of
\begin{equation*}\label{expression of T--11}T=-m(Q)+m(Q)Qm(Q),\end{equation*} we have
$T\geq 0$ by Lemma~\ref{thm:mQm geq m}, and
\begin{align*}
  T^2+T & =-m(Q)Qm(Q)+m(Q) Q\cdot m(Q)\cdot Q m(Q) \\
   & =-m(Q)Q\big[I-m(Q)\big]Qm(Q).
\end{align*}
Let
$$R=m(Q)Q\big[I-m(Q)\big].$$
Since $\big(m(Q),Q\big)$ is a quasi-projection pair, we have
$$R^*=\big[I-m(Q)\big]Q^*m(Q)=-\big[I-m(Q)\big]Qm(Q).$$
Therefore
$$RR^*=-m(Q)Q\big[I-m(Q)\big]Qm(Q)=T^2+T,$$
which is combined with the positivity of $T$ to get
\begin{equation*}\label{equ:norm relation of R and T}
\|R\|^2=\|T\|^2+\|T\|.
\end{equation*}
By Theorem~\ref{thm:matched projection for I-Q}, the operators $S$ and $R^*$ can be expressed alternatively as
\begin{align*}&S=-m(I-Q)\left[I-(I-Q)\right]m(I-Q),\\
&R^*=m(I-Q)(I-Q)\left[I-m(I-Q)\right].
\end{align*}
So if we replace $Q$ with $I-Q$, then from the  argument as above it can be concluded that $S\ge 0$ and $R^*R=S^2+S$. Consequently, it also holds that
$$\|R\|^2=\|S\|^2+\|S\|.$$
Hence,
$$\|T\|^2+\|T\|=\|S\|^2+\|S\|,$$
which happens only if $\|T\|=\|S\|$.
\end{proof}

\begin{lemma}\label{lem:norm relation about mq1-mq2 -02} Let $Q_1,Q_2\in \mathcal{L}(H)$ be idempotents such that  $$\alpha:= \left\|\big[I-m(Q_1)\big]Q_1\big[I-m(Q_1)\big]\right\|<1.$$
Then
\begin{equation}\label{equ:norm relation of mq12 and q12}
  \big\|m(Q_1)-m(Q_2)\big\|\leq \frac{\|Q_1-Q_2\|}{1-\alpha}.
\end{equation}
\end{lemma}
\begin{proof} For each idempotent $Q\in \mathcal{L}(H)$, let
$$Y=m(Q)Qm(Q)+\big[I-m(Q)\big]Q\big[I-m(Q)\big].$$
Since  $\big[I-m(Q)\big]\big[m(Q)Qm(Q)\big]^\dag=0$, we have
$$m(Q)Qm(Q)\big[m(Q)Qm(Q)\big]^\dag=Y \big[m(Q)Qm(Q)\big]^\dag.$$
As $\big(m(Q),Q\big)$ is a quasi-projection pair, $Y$ can be rewritten as
\begin{equation*}
   Y=\frac{1}{2}m(Q)(Q+Q^*)m(Q)+\frac{1}{2}\big[I-m(Q)\big](Q+Q^*)\big[I-m(Q)\big].\end{equation*}
By \eqref{equ:m qstar} and \eqref{equ:m times q}, we see that
$m(Q)(Q+Q^*)$ is self-adjoint, hence
\begin{equation}\label{equ:commutaty of mq and qqstar}
  m(Q)(Q+Q^*)=(Q+Q^*)m(Q).
\end{equation}
It follows that
\begin{align*}
   Y=&\frac{1}{2}m(Q)(Q+Q^*)+\frac{1}{2}\big[I-m(Q)\big](Q+Q^*)=\frac{1}{2}(Q+Q^*).
\end{align*}
Therefore,
\begin{equation*}m(Q)=m(Q)Qm(Q)\big[m(Q)Qm(Q)\big]^\dag=\frac{1}{2}(Q+Q^*)\big[m(Q)Qm(Q)\big]^\dag.\end{equation*}

Now, let $V_1$ and $V_2$ be defined as in the proof of Theorem~\ref{thm:norm relation about mq1-mq2 -01}. The above argument shows that
$$m(Q_2)=V_2V_2^\dag=\frac{1}{2}(Q_2+Q_2^*)V_2^\dag,\quad V_1+Z=\frac{1}{2}(Q_1+Q_1^*),$$
where $Z$ is defined by  $$Z=\big[I-m(Q_1)\big]Q_1\big[I-m(Q_1)\big]$$  such that $\alpha:=\|Z\|<1$, as is assumed.
In view of $$\big[I-m(Q_1)\big]V_1=\big[I-m(Q_2)\big]V_2^\dag=0,$$ we have
\begin{align*}0=&\big[I-m(Q_1)\big]\Big[V_1+Z\big[I-m(Q_2)\big]\Big]V_2^\dag\\
=&\big[I-m(Q_1)\big]\Big[\frac12(Q_1+Q_1^*)-W\Big]V_2^\dag,
\end{align*}
in which
$$W=Zm(Q_2)=Z\big[I-m(Q_1)\big]m(Q_2),$$
since $Z=Z\big[I-m(Q_1)\big]$. Hence,
\begin{align*}
  \big[I-m(Q_1)\big]m(Q_2)=\big[I-m(Q_1)\big]\Big[\frac{Q_2+Q_2^*}{2}-\frac{Q_1+Q_1^*}{2}+W\Big]V_2^\dag,
\end{align*}
which gives
\begin{align*}
  \big\|\big[I-m(Q_1)\big]m(Q_2)\big\| & \leq \left\|\frac{Q_2+Q_2^*}{2}-\frac{Q_1+Q_1^*}{2}+W\right\| \leq \|Q_1-Q_2\|+\|W\|\\
   & \leq \|Q_1-Q_2\|+\alpha \cdot \big\|\big[I-m(Q_1)\big]m(Q_2)\big\|.
\end{align*}
As a result,
$$\big\|\big[I-m(Q_1)\big]m(Q_2)\big\|\leq \frac{\|Q_1-Q_2\|}{1-\alpha}.$$
In virtue of Theorem~\ref{thm:matched projection for I-Q} and Lemma~\ref{lem:norm equation about matched projection}, we may replace
$Q_i$ with $I-Q_i$ for $i=1,2$ to
obtain
$$\big\|m(Q_1)\big[I-m(Q_2)\big]\big\|\leq \frac{\|Q_1-Q_2\|}{1-\alpha}.$$
The desired conclusion is immediate from Krein-Krasnoselskii-Milman equality.
\end{proof}

Now, we are in the position to give a partial answer to Problem~\ref{problem m Q1 Q2}.
\begin{theorem}\label{thm:conjecture inequality-special}For every projection $P\in \mathcal{L}(H)$ and idempotent $Q\in\mathcal{L}(H)$, we have
\begin{equation}\label{conjecture inequality-special} \|P-m(Q)\|\le \|P-Q\|.\end{equation}
\end{theorem}
\begin{proof}Let
$Q_1=P$ and $Q_2=Q$. In this case $m(Q_1)=Q_1$, so  the number $\alpha$ defined in Lemma~\ref{lem:norm relation about mq1-mq2 -02} turns out to be zero. The desired conclusion is immediate from
\eqref{equ:norm relation of mq12 and q12}.
\end{proof}
\begin{remark}\label{rem:conjecture inequality-special}
A simple use of Theorem~\ref{thm:conjecture inequality-special} yields
$$\big\|m(Q_1)-m(Q_2)\big\|\leq \min\big\{\|m(Q_1)-Q_2\|, \|Q_1-m(Q_2)\|\big\}$$
for every idempotents $Q_1,Q_2\in \mathcal{L}(H)$.
\end{remark}

We end this subsection by providing another partial answer to Problem~\ref{problem m Q1 Q2}.
\begin{theorem}\label{thm:norm up of mq1-mq2 accosiated with two projections} Let $P_1,P_2\in \mathcal{L}(H)$ be projections such that $\|P_1P_2\|< 1$. Then
$\|m(Q_1)-m(Q_2)\|\le \|Q_1-Q_2\|$, where $Q_1$ and $Q_2$ are idempotents defined by
\begin{equation}\label{equ:new idempotent Q1 and Q2}
Q_1=(I-P_1P_2)^{-1}P_1(I-P_2),\quad Q_2=(I-P_2P_1)^{-1}P_2(I-P_1).
\end{equation}
\end{theorem}
\begin{proof} In view of $P_1(I-P_2)=P_1(I-P_1P_2)$ and $P_2(I-P_1)=P_2(I-P_2P_1)$, it is easily seen that $Q_1$ and $Q_2$ defined by \eqref{equ:new idempotent Q1 and Q2}
are idempotents. If $P_1=0$ or $P_2=0$, then both of $Q_1$ and $Q_2$ are projections, so in this case the assertion is trivially satisfied.

Assume now that both of $P_1$ and $P_2$ are non-zero projections. A simple use of Krein-Krasnoselskii-Milman equality (see \eqref{equ:KKM equality}) shows that
$\|m(Q_1)-m(Q_2)\|\le 1$, so it needs only to prove that $\|Q_1-Q_2\|\ge 1.$ For this, we put
$S=S_1+S_2$,
where
$$S_1=P_1(I-P_2)P_1,\quad S_2=(I-P_1)P_2(I-P_1).$$
Evidently,
$$S_1\ge 0,\quad S_2\ge 0,\quad \mathcal{R}(S)=\mathcal{R}(S_1)+\mathcal{R}(S_2),\quad \overline{\mathcal{R}(S)}=\overline{\mathcal{R}(S_1)}+\overline{\mathcal{R}(S_2)}.$$
As $I-P_2P_1$ is invertible in $\mathcal{L}(H)$, we have
\begin{align*}\mathcal{R}(S_1)=\mathcal{R}\big[P_1(I-P_2P_1)\big]=\mathcal{R}(P_1),\end{align*}
so $\mathcal{R}(S_1)$ is closed in $H$. Similarly,
$\mathcal{R}\big[P_2(I-P_1)\big]$ is closed in $H$, which  leads by \cite[Lemma~5.7]{Luo-Moslehian-Xu} to the closedness of $\mathcal{R}(S_2)$ in $H$, since
$S_2=\big[P_2(I-P_1)\big]^*\big[P_2(I-P_1)\big]$. Therefore, $\mathcal{R}(S)$ is closed in $H$ and thus by \cite[Theorem~2.2]{Xu-sheng} $S$ is Moore-Penrose invertible. Let $S^\dag$ denote the Moore-Penrose inverse of $S$.
Then $S^\dag\ge 0$ and
$$SS^\dag=S^\dag S=P_{\mathcal{R}(S_1)+\mathcal{R}(S_2)}=P_{\mathcal{R}(P_1)+\mathcal{R}(S_2)},$$
so $SS^\dag P_1=S^\dag SP_1=P_1$, which gives $P_1SS^\dag =P_1S^\dag S=P_1$ by taking $*$-operation.
From the definitions of $S,S_1$ and $S_2$, it is clear that
$$P_1S=SP_1=S_1=(I-P_1P_2)P_1,$$ so $(I-P_1P_2)^{-1}P_1S=P_1$, and thus
\begin{align}\label{(i-PQ)inv P=PSdag} & (I-P_1P_2)^{-1}P_1=(I-P_1P_2)^{-1}P_1SS^\dag =P_1S^\dag,\\
  & S^\dag P_1=S^\dag P_1SS^\dag=S^\dag S P_1S^\dag=P_1S^\dag.\nonumber
 \end{align}
It follows from \eqref{equ:new idempotent Q1 and Q2} and \eqref{(i-PQ)inv P=PSdag} that
\begin{align*}
  Q_1Q_1^* & =(I-P_1P_2)^{-1}P_1(I-P_2) P_1 (I-P_2P_1)^{-1} \\
  & =(I-P_1P_2)^{-1}P_1(I-P_2P_1)(I-P_2P_1)^{-1}\\
  &=(I-P_1P_2)^{-1}P_1=P_1S^\dag.
\end{align*}
Also, from the definitions of $S,S_1$ and $S_2$ we have
\begin{align*}S=P_1+P_2-P_1P_2-P_2P_1=P_2(I-P_1)P_2+(I-P_2)P_1(I-P_2),
\end{align*}
so we may exchange $P_1$ with $P_2$ to obtain
\begin{equation}\label{equ: the property of Q2}
   S^\dag P_2=P_2S^\dag,\quad (I-P_2P_1)^{-1}P_2=P_2S^\dag,\quad Q_2Q_2^*=P_2S^\dag.
\end{equation}
Utilizing \eqref{equ:new idempotent Q1 and Q2}   yields
\begin{align*}Q_1Q_2^*=&(I-P_1P_2)^{-1}P_1(I-P_2)(I-P_1)P_2(I-P_1P_2)^{-1}\\
=&-(I-P_1P_2)^{-1}P_1P_2(I-P_1)P_2(I-P_1P_2)^{-1}\\
=&-(I-P_1P_2)^{-1}P_1P_2(I-P_1P_2)(I-P_1P_2)^{-1}\\
=&-(I-P_1P_2)^{-1}P_1P_2.\end{align*}
So by \eqref{(i-PQ)inv P=PSdag} and the first equation in \eqref{equ: the property of Q2}, we have
$$Q_1Q_2^*=-P_1S^\dag P_2=-P_1P_2 S^\dag.$$
Similarly, $Q_2Q_1^*=-P_2P_1 S^\dag$. Consequently,
\begin{align*}(Q_1-Q_2)(Q_1-Q_2)^*=&Q_1Q_1^*-Q_1Q_2^*-Q_2Q_1^*+Q_2Q_2^*\\
=&P_1S^\dag+P_1P_2 S^\dag+P_2P_1 S^\dag+P_2S^\dag=TS^\dag,
\end{align*}
where $$T=P_1+P_1P_2+P_2P_1+P_2.$$
Employing $P_1=P_1SS^\dag$ and $P_1S=P_1-P_1P_2P_1$, we see that
\begin{align*}
  (Q_1-Q_2)(Q_1-Q_2)^*-P_1 & =TS^\dag -P_1SS^\dag=(T-P_1+P_1P_2P_1)S^\dag\\
   & =(I+P_1)P_2(I+P_1)S^\dag\\
   &=\big(S^\dag\big)^\frac12(I+P_1)P_2(I+P_1)\big(S^\dag\big)^\frac12 \geq 0.
\end{align*}
Therefore,
$\|Q_1-Q_2\|\geq \sqrt{\|P_1\|}= 1.$
\end{proof}

\subsection{Investigations  of  $\|m(Q)-Q\|\le \|P-Q\|$}

Given an arbitrary idempotent $Q\in M_2(\mathbb{C})$, it has already been shown in Section~\ref{sec:quasi-projection pairs} that
$\|m(Q)-Q\|\le \|P-Q\|$ for every projection $P\in M_2(\mathbb{C})$. This leads us to put forward the following problem in the general setting of Hilbert $C^*$-module operators.

\begin{problem}\label{problem m Q Q P} Under what conditions it is true that $\|m(Q)-Q\|\leq \|P-Q\|$ for a projection $P$  and an idempotent $Q$ in  $\mathcal{L}(H)$?\end{problem}

To get a partial answer to the above problem, we begin with the special case that $P=P_{\mathcal{R}(Q)}$. For this, we provide a lemma as follows.

\begin{lemma}\label{lem:significance} For every idempotent $Q\in \mathcal{L}(H)$, we have
\begin{align*}
   &(Q^*-Q)(Q^*-Q)^*=4D^2+4D,\\
   &\big[m(Q)-Q\big] \big[m(Q)-Q\big]^*+\big[m(Q)-Q\big]^*\big[m(Q)-Q\big]=4D^2+2D,
\end{align*}
where $m(Q)$ is the matched projection of $Q$ and
$$D=-m(Q)(I-Q)m(Q)-\big[I-m(Q)\big]Q\big[I-m(Q)\big].$$
\end{lemma}
\begin{proof}  A simple calculation shows that
\begin{align*}D=&m(Q)Q+Qm(Q)-Q-m(Q)\\
=&\frac{1}{2}\big[2m(Q)-I\big]Q+\frac{1}{2}Q\big[2m(Q)-I\big]-m(Q).\end{align*}
So, it can be derived from \eqref{qpp for m Q}  that $D=D_1=D_2$, where
$$D_1=\frac{1}{2}(Q+Q^*)\big[2m(Q)-I]-m(Q), \quad D_2=\frac{1}{2}\big[2m(Q)-I](Q+Q^*)-m(Q).$$
In virtue of \eqref{equ:commutaty of mq and qqstar}, we have
\begin{align*}
  D^2 & =D_1D_2=\frac{1}{4}(Q+Q^*)^2-\frac{1}{2}(Q+Q^*)m(Q)-\frac{1}{2}m(Q)(Q+Q^*)+m(Q) \\
   & =\frac{1}{4}(Q+Q^*)^2-m(Q)(Q+Q^*)+m(Q).
\end{align*}
Therefore,
\begin{align*}
  4D^2+4D & =4D^2+4D_2=(Q+Q^*)^2-2(Q+Q^*) \\
   & =QQ^*+Q^*Q-Q-Q^*=(Q^*-Q)(Q^*-Q)^*.
\end{align*}
To verify the second assertion, we put
\begin{equation}\label{equ:auxiliary positive operators X and Y}
X=\big[m(Q)-Q\big]\big[m(Q)-Q\big]^*,\quad Y=\big[m(Q)-Q\big]^*\big[m(Q)-Q\big].
\end{equation}
Once again, by  \eqref{equ:commutaty of mq and qqstar} we have
\begin{align*}
  X+Y & =2m(Q)-(Q+Q^*)m(Q)-m(Q)(Q+Q^*)+QQ^*+Q^*Q \\
   & =2m(Q)-2m(Q)(Q+Q^*)+QQ^*+Q^*Q\\
   &=2m(Q)-\big[2m(Q)+I](Q+Q^*)+(Q+Q^*)^2\\
   &=4D^2+2D_2=4D^2+2D.
\end{align*}
This completes the proof.
\end{proof}

Our first partial answer to Problem~\ref{problem m Q Q P} is as follows.

\begin{theorem}\label{thm:significance} For every idempotent $Q\in \mathcal{L}(H)$, we have
\begin{equation}\label{equ:calimed sharpness}\frac12 \left\|P_{\mathcal{R}(Q)}-Q\right\|\le \left\|m(Q)-Q\right\|\le\left\|P_{\mathcal{R}(Q)}-Q\right\|.
\end{equation}
Moreover, each equality above occurs if and only if $Q$ is a projection.
\end{theorem}
\begin{proof} It is well-known (see e.g.\,\cite[Proposition~4.3]{Koliha}) that $\|P_{\mathcal{R}(Q)}-Q\|=\|Q^*-Q\|$.
So, it needs only to verify that
\begin{equation}\label{lub 1 2}\|m(Q)-Q\|\le \|Q^*-Q\|\le 2\|m(Q)-Q\|.\end{equation}
Following the notations as in the proof of Lemma~\ref{lem:significance}, we have
\begin{equation}\label{equ:relation of four positive operators}
X+Y=4D^2+2D,\quad (Q^*-Q)(Q^*-Q)^*=4D^2+4D,
\end{equation}
where $X$ and $Y$ are given by \eqref{equ:auxiliary positive operators X and Y} and
$D=T+S$, in which $T$ and $S$ are defined by \eqref{equ:auxiliary positive operators}. From the proof of Lemma~\ref{lem:norm equation about matched projection}, we have $T\geq 0$ and $S\geq0$, hence  $D\geq 0$. The positivity of $D$ together with \eqref{equ:auxiliary positive operators X and Y} and \eqref{equ:relation of four positive operators} yields
\begin{align}\label{norm X and Y}&\|X\|=\|Y\|=\|m(Q)-Q\|^2,\\
\label{norm related eith D}&\|X+Y\|=4\|D\|^2+2\|D\|,\quad \|Q^*-Q\|^2=4\|D\|^2+4\|D\|,\\
&0\le X\leq X+Y\leq (Q^*-Q)(Q^*-Q)^* \leq 2(X+Y).\nonumber
\end{align}
Hence
\begin{equation*}
\|X\|\leq \|X+Y\|\leq \|Q^*-Q\|^2\leq 2\|X+Y\|\leq 2(\|X\|+\|Y\|),
\end{equation*}
which is combined with \eqref{norm X and Y} and \eqref{norm related eith D} to get
\begin{equation}\label{equ:norm relation of four positive operators}\|m(Q)-Q\|^2\le \|X+Y\|\leq \|Q^*-Q\|^2\le 2\|X+Y\| \leq 4\|m(Q)-Q\|^2.\end{equation}
This shows the validity of \eqref{lub 1 2}.

Clearly, all the norms in \eqref{lub 1 2} are equal to zero whenever $Q$ is a projection. Suppose that $\|m(Q)-Q\|=\|Q^*-Q\|$. In this case, by \eqref{equ:norm relation of four positive operators}
we have $\|X+Y\|=\|Q^*-Q\|^2$, so a simple use of \eqref{norm related eith D} will lead to $D=0$, which in turn gives $Q^*=Q$ by \eqref{norm related eith D}. Similarly, if $\|Q^*-Q\|= 2\|m(Q)-Q\|$, then it can also be deduced from \eqref{equ:norm relation of four positive operators}
and \eqref{norm related eith D} that $Q^*=Q$.
\end{proof}

\begin{remark}\label{rem:optical of two numbers}By Theorem~\ref{thm:significance}, there exist constants $\alpha>0$ and $\beta>0$ such that
\begin{equation}\label{equ:calimed sharpness max min}\alpha\left\|P_{\mathcal{R}(Q)}-Q\right\|\le \left\|m(Q)-Q\right\|\le \beta \left\|P_{\mathcal{R}(Q)}-Q\right\|
\end{equation}
for all idempotent $Q\in\mathcal{L}(H)$.
For each $a>0$, let $Q$ be given as \eqref{idempotent in M two}. In this case
$m(Q)=P_0$, which is formulated by \eqref{equ:problem P0} with $b=\sqrt{a^2+1}$ therein.
 It is clear that
$\|P_{\mathcal{R}(Q)}-Q\|=a$.  Direct calculations show that
$$\big[m(Q)-Q\big]\big[m(Q)-Q\big]^*=\frac{b-1}{2b}\left(
                               \begin{array}{cc}
                                 2b^2-1 & -a \\
                                 -a & 1 \\
                               \end{array}
                             \right),$$
and the two eigenvalues of $\big[m(Q)-Q\big]\big[m(Q)-Q\big]^*$ are $\lambda_1=\frac{1}{2}(b-1)(b-a)$ and $\lambda_2=\frac{1}{2}(b-1)(b+a)$. Hence
$$\|m(Q)-Q\|^2=\left\|\big[m(Q)-Q\big]\big[m(Q)-Q\big]^*\right\|=\lambda_2.$$
Therefore
$$ \frac{\|m(Q)-Q\|^2}{\|P_{\mathcal{R}(Q)}-Q\|^2}=\frac{\lambda_2}{a^2}=\frac{\sqrt{a^2+1}+a}{2\big(\sqrt{a^2+1}+1\big)},$$
which gives
$$\lim_{a\to 0} \frac{\|m(Q)-Q\|}{\|P_{\mathcal{R}(Q)}-Q\|}=\frac{1}{2},\quad \lim_{a\to \infty} \frac{\|m(Q)-Q\|}{\|P_{\mathcal{R}(Q)}-Q\|}=1.$$
The limitations above together with \eqref{equ:calimed sharpness} indicate that $\frac12$ and $1$ are the optimal numbers of $\alpha$ and $\beta$ in
\eqref{equ:calimed sharpness max min},  respectively.
\end{remark}

If we take $P=0$, $P=I$ and $P=P_{\mathcal{N}(Q)}$ respectively for an arbitrary idempotent $Q\in\mathcal{L}(H)$, then we have the following corollary, which also gives a partial answer to Problem~\ref{problem m Q Q P}.
\begin{corollary}\label{cor:partial answer of pro.2} For every idempotent $Q\in\mathcal{L}(H)$, we have
$$\|m(Q)-Q\|\le  \|Q\|=\|I-Q\|\le \|P_{\mathcal{N}(Q)}-Q\|.$$
\end{corollary}
\begin{proof}
 Let $Q$ be an arbitrary idempotent on $H$. We claim that
\begin{equation*}\label{equ:relationship norm Q}
  \|P_{\mathcal{N}(Q)}-Q\|\ge \|Q\|= \|I-Q\|\ge \|P_{\mathcal{R}(Q)}-Q\|\ge \|m(Q)-Q\|.
\end{equation*}
Indeed, the equation $\|Q\|= \|I-Q\|$ is known (see e.\,g.\,\cite[Proposition~4.2]{Koliha}), and the last inequality has already been proved (see \eqref{equ:calimed sharpness}). To prove the rest inequalities, we notice that
$$QP_{\mathcal{N}(Q)}=0=P_{\mathcal{N}(Q)}Q^*,\quad QP_{\mathcal{R}(Q)}=P_{\mathcal{R}(Q)}=P_{\mathcal{R}(Q)}Q^*,$$
so
\begin{align*}
   & \big(P_{\mathcal{N}(Q)}-Q\big)\big(P_{\mathcal{N}(Q)}-Q\big)^*=P_{\mathcal{N}(Q)}+QQ^*, \\
   & \big(P_{\mathcal{R}(Q)}-Q\big)\big(P_{\mathcal{R}(Q)}-Q\big)^*=-P_{\mathcal{R}(Q)}+QQ^*.
\end{align*}
Hence
$$\big(P_{\mathcal{N}(Q)}-Q\big)\big(P_{\mathcal{N}(Q)}-Q\big)^*\ge QQ^* \ge \big(P_{\mathcal{R}(Q)}-Q\big)\big(P_{\mathcal{R}(Q)}-Q\big)^*,$$
which means that
$$\|P_{\mathcal{N}(Q)}-Q\|\ge \|Q\|\ge \|P_{\mathcal{R}(Q)}-Q\|.\qedhere$$
\end{proof}

To provide another partial answer to Problem~\ref{problem m Q Q P}, we need a technical lemma as follows.

\begin{lemma}\label{thm:mq-q times mq-qstar is min} For every quasi-projection pair $(P,Q)$ on $H$, we have
\begin{align}\label{equ:mq-q times mq-qstar is min}&
  (P-Q)(P-Q)^*\geq \big[m(Q)-Q\big]\big[m(Q)-Q\big]^*,\\
 \label{equ:mq-q times mq-qstar is min+} &(P-Q)^*(P-Q)\geq \big[m(Q)-Q\big]^*\big[m(Q)-Q\big].
\end{align}
Moreover, each equality occurs if and only if $P=m(Q)$.
\end{lemma}
\begin{proof} It needs to show that the operator $W\in\mathcal{L}(H)$ defined by
$$W=(P-Q)(P-Q)^*- \big[m(Q)-Q\big]\big[m(Q)-Q\big]^*$$
is positive. Direct computations yield
$$W=P-m(Q)-Q\big[P-m(Q)\big]-\big[P-m(Q)\big]Q^*.$$
Since  both $(P,Q)$ and  $\big(m(Q),Q\big)$ are quasi-projection pairs, we may use Theorem~\ref{thm:short description of qpp} together with the decomposition
 $$P-m(Q)=\frac{1}{2}\left[(2P-I)-\big(2m(Q)-I\big)\right]$$
 to get
\begin{align}\label{semi-commutative-002}
  Q\big[P-m(Q)\big] & = \frac{1}{2}\left[(2P-I)-\big(2m(Q)-I\big)\right]Q^*=\big[P-m(Q)\big]Q^*.
\end{align}
Therefore, the operator $W$ defined as above can be simplified as
\begin{equation}\label{equ:preliminary calculation of W}
  W=-\big[P-m(Q)\big](2Q^*-I).
\end{equation}
By Corollary~\ref{cor:commutativity of P and mq} we have $Pm(Q)=m(Q)P$, which yields
\begin{equation}\label{equ:the relationship of p-mq and 2m-1}
  \big[P-m(Q)\big]\cdot \big[2m(Q)-I\big]\cdot\big[P-m(Q)\big]=-\big[P-m(Q)\big].
\end{equation}
Moreover, since $m(I-Q)=I-m(Q)$ (see Theorem~\ref{thm:matched projection for I-Q}), it can be derived directly from  \eqref{equ:relation of q and 1q1} that
$$\big[2m(Q)-I\big](Q-I)=|I-Q|$$
by replacing $Q$ therein with $I-Q$.
The equation above together with \eqref{equ:relation of q and 1q1} yields
\begin{equation}\label{equ:positive operator 2mq-1 times 2q-1}
 \big[2m(Q)-I\big](2Q-I)=|Q|+|I-Q|.
\end{equation}
A combination of  \eqref{equ:positive operator 2mq-1 times 2q-1}, \eqref{semi-commutative-002},  \eqref{equ:the relationship of p-mq and 2m-1} and \eqref{equ:preliminary calculation of W} gives
\begin{align}\label{equ:preliminary calculation of W 02}
 \big[P-m(Q)\big]\cdot (|Q|+|I-Q|)\cdot \big[P-m(Q)\big]=W.
\end{align}
This shows the positivity of $W$.

Clearly, the operator $\big[2m(Q)-I\big](2Q-I)$ is invertible in $\mathcal{L}(H)$,
so it can be deduced from \eqref{equ:positive operator 2mq-1 times 2q-1} that $|Q|+|I-Q|$ is positive definite.
Hence, by \eqref{equ:preliminary calculation of W 02} we see that
$$W=0 \Longleftrightarrow (|Q|+|I-Q|)\cdot \big[P-m(Q)\big]=0 \Longleftrightarrow P-m(Q)=0. $$
Therefore, the equality in \eqref{equ:mq-q times mq-qstar is min} occurs if and only if $P=m(Q)$.

In virtue of Theorems~\ref{thm:four equivalences} and \ref{thm:same matched projection for Q and Q-sta}, \eqref{equ:mq-q times mq-qstar is min+}
is immediate from \eqref{equ:mq-q times mq-qstar is min} by replacing the pair $(P,Q)$ with $(P,Q^*)$.
\end{proof}

Our next partial answer to Problem~\ref{problem m Q Q P}  reads as follows.

\begin{theorem}\label{thm:m Q is partially best} For every projection $P\in \mathcal{L}(H)$ and idempotent $Q\in \mathcal{L}(H)$, we have
\begin{equation}\label{two times not so good}
  \|m(Q)-Q\|\leq 2\|P-Q\|.
\end{equation}
If furthermore $(P,Q)$ is a quasi-projection pair, then $\|m(Q)-Q\|\leq \|P-Q\|$.
\end{theorem}
\begin{proof}It follows from \eqref{conjecture inequality-special} that
$$\|m(Q)-Q\|\leq \|m(Q)-P\|+\|P-Q\|\leq 2\|P-Q\|.$$
This shows the validity of \eqref{two times not so good}. In addition, if $(P,Q)$ is a quasi-projection pair, then a direct use of  \eqref{equ:mq-q times mq-qstar is min} yields $\|m(Q)-Q\|\leq \|P-Q\|$.
\end{proof}

\subsection{A formula for $\|Q-m(Q)\|$ and its applications}
It is interesting to provide a formula for
$\|m(Q)-Q\|$ in the case that $Q$ is a  general idempotent on a Hilbert $C^*$-module. For this, we need a lemma as follows.

\begin{lemma}\label{lem:pre case p-mq le 1} For each idempotent $Q\in \mathcal{L}(H)$, let
\begin{equation*}\label{equ:positive defined operator v}
V=\frac{1}{2}\big[|Q|+|I-Q|+I\big].
\end{equation*}
Then $V$ is a positive definite operator which satisfies
 $Q=V^{-1}m(Q)V$ and
$$(I-V)^2=\big[m(Q)-Q\big]^*\big[m(Q)-Q\big].$$
\end{lemma}
\begin{proof} Obviously, the operator $V$ defined as above is positive definite.  It can be deduced from \eqref{equ:positive operator 2mq-1 times 2q-1} that
$$V=2m(Q)Q-m(Q)-Q+I,$$
which clearly leads to
$$VQ=m(Q)Q=m(Q)V.$$
Therefore, $Q=V^{-1}m(Q)V$. Furthermore, the latter expression of $V$ together with \eqref{qpp for m Q} yields  $I-V=A=B$, where
$$A=-\big[2m(Q)-I\big]Q+m(Q),\quad B=-Q^*\big[2m(Q)-I\big]+m(Q).$$
Consequently,
\begin{align*}(I-V)^2=&BA=m(Q)-Q^*m(Q)-m(Q)Q+Q^*Q\\
=&\big[m(Q)-Q\big]^*\big[m(Q)-Q\big].\qedhere\end{align*}
\end{proof}

We are now in the position to provide a formula for $\|m(Q)-Q\|$.

\begin{theorem}\label{thm:norm of mq-q} For each idempotent $Q\in \mathcal{L}(H)$, we have
\begin{equation}\label{norm of q minus m q}\|m(Q)-Q\|=\frac{1}{2}\left[\|Q\|-1+\sqrt{\|Q\|^2-1}\right].\end{equation}
\end{theorem}
\begin{proof}Apparently, the conclusion is true whenever $Q$ is a projection. So, in what follows we assume that $Q$ is a non-projection idempotent.
In view of Theorem~\ref{thm:invance of matched projection}, we have
 $\|m(Q)-Q\|=\|m\big(\pi(Q)\big)-\pi(Q)\|$ and $\|Q\|=\|\pi(Q)\|$ for every Hilbert space $X$ that ensures a faithful unital $C^*$-morphism   $\pi:\mathcal{L}(H)\to\mathbb{B}(X)$. So, without loss of generality, we may as well assume that $H$ is a Hilbert space. In this case, every bounded linear operator considered has the polar decomposition.

 By Lemma~\ref{lem:pre case p-mq le 1}, we have
\begin{equation}\label{preformula for m q minus q-01}\|m(Q)-Q\|=\|V-I\|=\frac{1}{2}\big\||Q|+|I-Q|-I\big\|.\end{equation}
Following the notations as in the proof of Theorem~\ref{existence thm of quasi-projection pair}, we see that
$\widetilde{Q}$ is defined by  \eqref{defn of widetile Q-1}, in which $A$  is given by  \eqref{equ:1st defn of AB}.
It is clear that
\begin{align*}\|Q\|=\|\widetilde{Q}\|=\sqrt{\|\widetilde{Q}\widetilde{Q}^*}\|=\sqrt{\|I_1+AA^*\|}=\sqrt{1+\|AA^*\|}=\sqrt{1+\|A\|^2}.
\end{align*}
Hence,
\begin{equation}\label{norm A wrt Q}\|A\|=\sqrt{\|Q\|^2-1}.\end{equation}
Direct computations yield
\begin{align*}&
|\widetilde{Q}|=\left(
                               \begin{array}{cc}
                                 (I_1+AA^*)^{-\frac12} & (I_1+AA^*)^{-\frac12}A\\
                                 A^*(I_1+AA^*)^{-\frac12} & A^*(I_1+AA^*)^{-\frac12}A\\
                               \end{array}
                             \right),\\
&|\widetilde{I}-\widetilde{Q}|=\left(
                                 \begin{array}{cc}
                                   0 & 0 \\
                                   0 & (I_2+A^*A)^\frac12 \\
                                 \end{array}
                               \right).
\end{align*}
Since $Ap(A^*A)=p(AA^*)A$ for every polynomial $p$, we have
$$Af(A^*A)=f(AA^*)A$$ whenever $f$ is a function which is continuous on $[0,\|A\|^2]$.  Consequently,
\begin{align*}(I_2+A^*A)^\frac12 - I_2=&A^*A\left[(I_2+A^*A)^\frac12+I_2\right]^{-1}\\
=&A^*\left[(I_1+AA^*)^\frac12+I_1\right]^{-1}A.
\end{align*}
Let
$$T=|\widetilde{Q}|+|\widetilde{I}-\widetilde{Q}|-\widetilde{I},$$
and let   $A=V_A|A|$ be the polar decomposition of $A$ \cite[Lemma~3.6]{Liu-Luo-Xu}. Then the above arguments show that
$$T=\left(
      \begin{array}{cc}
        I_1 & 0 \\
        0 & V_A^* \\
      \end{array}
    \right)S\left(
      \begin{array}{cc}
        I_1 & 0 \\
        0 & V_A \\
      \end{array}
    \right),$$
where $S=(S_{ij})_{1\le i,j\le 2}$, in which
\begin{align*}S_{11}=&(I_1+|A^*|^2)^{-\frac12}-I_1,\quad S_{12}=(I_1+|A^*|^2)^{-\frac12}|A^*|,\quad S_{21}=S_{12}^*,\\
      S_{22}=&|A^*|(I_1+|A^*|^2)^{-\frac12}|A^*|+|A^*|\left[(I_1+|A^*|^2)^\frac12+I_1\right]^{-1}|A^*|\\
      =&|A^*|^2\left[(I_1+|A^*|^2)^{-\frac12}+\left[(I_1+|A^*|^2)^\frac12+I_1\right]^{-1}\right].
      \end{align*}
Therefore,
$$\|T\|\le \left\|\left(
      \begin{array}{cc}
        I_1 & 0 \\
        0 & V_A^* \\
      \end{array}
    \right)\right\|\cdot \|S\|\cdot\left\|\left(
      \begin{array}{cc}
        I_1 & 0 \\
        0 & V_A \\
      \end{array}
    \right)\right\|\le \|S\|.$$
Conversely, due to $V_AV_A^*|A^*|=|A^*|$ and the above expressions of $S_{ij} (i,j=1,2)$, we have   $$\left(
      \begin{array}{cc}
        I_1 & 0 \\
        0 & V_A \\
      \end{array}
    \right)T\left(
      \begin{array}{cc}
        I_1 & 0 \\
        0 & V_A^* \\
      \end{array}
    \right)=\left(   \begin{array}{cc}
        I_1 & 0 \\
        0 & V_AV_A^* \\
      \end{array}
    \right)
    S\left(
      \begin{array}{cc}
        I_1 & 0 \\
        0 & V_AV_A^* \\
      \end{array}\right)=S,$$
which gives  $\|S\|\le \|T\|$.  So, according to \eqref{preformula for m q minus q-01} we have
\begin{equation}\label{preformula for m q minus q-02}\|m(Q)-Q\|=\frac12\|T\|=\frac12\|S\|.
\end{equation}
For each $t\in [0,+\infty)$, let
\begin{align*}S_t=\left(
                    \begin{array}{cc}
                      \frac{1}{\sqrt{1+t^2}}-1 & \frac{t} {\sqrt{1+t^2}}\\
                      \frac{t} {\sqrt{1+t^2} }& \frac{t^2}{\sqrt{1+t^2}}+\frac{t^2}{\sqrt{1+t^2}+1}\\
                    \end{array}
                  \right).
\end{align*}
Then for each $\lambda\in (-\infty,+\infty)$, we have
\begin{align*}\mbox{det}\left[\mbox{diag}(\lambda,\lambda)-S_t\right]=\lambda^2-2(\sqrt{1+t^2}-1)\lambda-2(\sqrt{1+t^2}-1).
\end{align*}
Therefore, $\lambda_i(t)(i=1,2)$ are the eigenvalues of $S_t$, where
\begin{align*}\lambda_1(t)=\sqrt{1+t^2}-1+t,\quad \lambda_2(t)=\sqrt{1+t^2}-1-t.
\end{align*}
Since $S_t$ is Hermitian and $t\ge 0$,  we have
$$\|S_t\|=\max\big\{|\lambda_1(t)|,|\lambda_2(t)|\big\}=\lambda_1(t).$$

Let $C(\sigma(|A^*|)$ be the set consisting of continuous functions on $\sigma(|A^*|)$ (the spectrum of $|A^*|$), and let $\mathfrak{B}$ be the unital commutative $C^*$-algebra generated by $I_1$ and $|A^*|$.
Observe that $\|A\|=\| \,|A^*|\,\|\in\sigma(|A^*|)$, and $\lambda_1(t)$ is monotonically increasing on $[0,+\infty)$, so it can be concluded that
\begin{equation}\label{preformula for m q minus q-03}\max_{t\in\sigma(|A^*|)} \lambda_1(t)=\lambda_1(\|A\|)=\sqrt{1+\|A\|^2}-1+\|A\|.\end{equation}
Since $S\in M_2(\mathfrak{B})\cong M_2\big[C\big(\sigma(|A^*|)\big)\big]\cong C\big(\sigma(|A^*|)\big)\otimes M_2(\mathbb{C})$,  and $M_2(\mathbb{C})$ is nuclear \cite[Proposition~T.5.20]{Wegge-Olsen}, we have
\begin{equation}\label{preformula for m q minus q-04}\|S\|=\max_{t\in\sigma(|A^*|)}\|S_t\|=\max_{t\in\sigma(|A^*|)} \lambda_1(t).
\end{equation}
Combining \eqref{preformula for m q minus q-02}, \eqref{preformula for m q minus q-04}, \eqref{preformula for m q minus q-03} and
\eqref{norm A wrt Q} yields the desired conclusion.
\end{proof}

\begin{remark}\label{rem: alternative proof} Suppose that $Q\in \mathcal{L}(H)$ is a non-projection idempotent.
From the proof of Theorem~\ref{thm:norm of mq-q}, we see that $\|Q\|=\sqrt{1+\|A\|^2}>1$ and
$$\|P_{\mathcal{R}(Q)}-Q\|=\left\|\left(
                             \begin{array}{cc}
                               I_1 & 0 \\
                               0 & 0 \\
                             \end{array}
                           \right)-\widetilde{Q}\right\|=\|A\|=\sqrt{\|Q\|^2-1}.$$
So, a simple use of \eqref{norm of q minus m q} will  lead to \eqref{equ:calimed sharpness}. This gives an alternative proof of Theorem~\ref{thm:significance}.                         \end{remark}

We provide a supplement to  Theorem~\ref{thm:m Q is partially best} as follows.

\begin{theorem}\label{thm:original NIE} Let $(P,Q)$ be a quasi-projection pair such that $\|P-Q\|<1$. Then $P=m(Q)$ and $\|Q\|<\frac53$.
\end{theorem}
\begin{proof} Suppose that $\|P-Q\|<1$. A direct application of \eqref{conjecture inequality-special} yields
$\big\|P-m(Q)\big\|<1$,
which in turn leads by Krein-Krasnoselskii-Milman equality  to
$$\|P_1\|<1, \quad \|P_2\|<1,$$ where
$$P_1=P\big[I-m(Q)\big],\quad P_2=(I-P)m(Q).$$
By Corollary~\ref{cor:commutativity of P and mq}  $P$ commutes with $m(Q)$,
so both of  $P_1$ and $P_2$ are projections.
Therefore, the latter two inequalities imply that $P_1=P_2=0$, which clearly gives $P=m(Q)$. As a result, we arrive at
$\|m(Q)-Q\|<1$. Utilizing \eqref{norm of q minus m q}, we see that
$$\frac{1}{2}\left[\|Q\|-1+\sqrt{\|Q\|^2-1}\right]<1,$$
which yields $\|Q\|<\frac53$, as desired.
\end{proof}

\vspace{2ex}
\noindent\textbf{Acknowledgement}

\vspace{2ex}

The authors thank the referee for helpful comments and suggestions.





\begin{thebibliography}{99}


\bibitem{Afriat} \textsc{S. Afriat}, \textit{Orthogonal and oblique projectors and the characteristics of pairs of vector spaces}, Proc. Cambridge Philos. Soc.,
\textbf{53}(1957), 800--816.



\bibitem{Ando02} \textsc{T. Ando}, \textit{Unbounded or bounded idempotent operators in Hilbert space}, Linear Algebra Appl., \textbf{438}(2013), no. 10,  3769--3775.





\bibitem{Aurel} \textsc{G. Aur$\acute{\mbox{e}}$l}, \textit{Projectors and projection methods}, Advances in Mathematics (Dordrecht),\textbf{6}. Kluwer Academic Publishers, Boston, MA 2004.






\bibitem{Blackadar}\textsc{B. Blackadar}, \textit{$K$-theory for operator algebras}, Mathematical Sciences Research Institute Publications, \textbf{5}, Springer-Verlag, New York 1986.



\bibitem{Blackadar-02}\textsc{B. Blackadar}, \textit{Projections in $C^*$-algebras, $C^*$-algebras: 1943--1993} (San Antonio, TX, 1993), 130--149, Contemp. Math., \textbf{167}, Amer. Math. Soc., Providence, RI, 1994.

\bibitem{Bottcher-Spitkovsky}\textsc{A.B\"ottcher, I. M. Spitkovsky}, \textit{A gentle guide to the basics of two projections theory}, Linear Algebra Appl., \textbf{432}(2010), no. 6,  1412--1459.




\bibitem{Corach}\textsc{G. Corach, H. Porta, L. Recht}, \textit{The geometry of spaces of projections in $C^*$-algebras}, Adv. Math., \textbf{101}(1993), 59--77.





\bibitem{Halmos} \textsc{P. Halmos}, \textit{Two subspaces}, Trans. Amer. Math. Soc., \textbf{144}(1969), 381--389.


\bibitem{Halpern} \textsc{H. Halpern, V. Kaftal, P. W. Ng, S. Zhang}, \textit{Finite sums of projections in von Neumann algebras}, Trans. Amer. Math. Soc., \textbf{365}(2013), no. 5,  2409--2445.







\bibitem{Kaftal} \textsc{V. Kaftal, P. W. Ng, S. Zhang}, \textit{Commutators and linear spans of projections in certain finite $C^*$-algebras}, J. Funct. Anal., \textbf{26}(2014), no. 4, 1883--1912.



\bibitem{Koliha} \textsc{J.J. Koliha}, \textit{Range projections of idempotents in $C^*$-algebras}, Demonstratio Math., \textbf{34}(2001), 91--103.




\bibitem{Lance}\textsc{E.C. Lance}, \textit{Hilbert $C^*$-modules--A toolkit for operator algebraists}, Cambridge University Press, Cambridge 1995.





\bibitem{Liu-Luo-Xu}\textsc{N. Liu, W. Luo, Q. Xu}, \textit{The polar decomposition for adjointable operators on Hilbert $C^*$-modules and centered operators}, Adv. Oper. Theory, \textbf{3}(2018), no. 4,  855--867.



\bibitem{Luo-Moslehian-Xu}\textsc{W. Luo, M.S. Moslehian, Q. Xu}, \textit{Halmos' two projections theorem for Hilbert $C^*$-module operators and the Friedrichs angle of two closed submodules}, Linear Algebra Appl., \textbf{577}(2019), 134--158.



\bibitem{MT} \textsc{V.M. Manuilov, E.V. Troitsky}, \textit{Hilbert $C^*$-modules}, Translated from the 2001 Russian original by the authors, Translations of Mathematical Monographs, \textbf{226}, American Mathematical Society, Providence, RI 2005.

\bibitem{Mesland}\textsc{B. Mesland, A. Rennie}, \textit{The Friedrichs angle and alternating projections in Hilbert $C^*$-modules}, J. Math. Anal. Appl., \textbf{516}(2022), Paper No. 126474, 19 pp.






\bibitem{Paschke}\textsc{W.L. Paschke}, \textit{Inner product modules over $B^*$-algebras}, Trans. Amer. Math. Soc., \textbf{182}(1973), 443--468.



\bibitem{Paulsen} \textsc{V.I. Paulsen}, \textit{Operator algebras of idempotents}, J. Funct. Anal., \textbf{181}(2001), no. 2, 209--226.


\bibitem{Pedersen02}\textsc{G.K. Pedersen}, \textit{Measure theory for $C^*$-algebras, II}, Math. Scand., \textbf{22}(1968), 63--74.


\bibitem{Raeburn}\textsc{I. Raeburn, A.M. Sinclair}, \textit{The $C^*$-algebra generated by two projections}, Math. Scand., \textbf{65}(1989), no. 2,  278--290.




\bibitem{Roch}\textsc{S. Roch, B. Silbermann}, \textit{Algebras generated by idempotents and the symbol calculus for singular integral operators}, Integral
Equations Operator Theory, \textbf{11} (1988), 385--419.

\bibitem{RLL}\textsc{M.R{\o}rdam, F. Larsen, N. Laustsen}, \textit{An introduction to $K$-theory for $C^*$-algebras}, London Mathematical Society Student Texts, \textbf{49}, Cambridge University Press, Cambridge 2000.




\bibitem{Wegge-Olsen}\textsc{N.E. Wegge-Olsen}, \textit{$K$-theory and $C^*$-algebras: a friendly approach},
Oxford University Press, New York 1993.



\bibitem{Xu-sheng}\textsc{Q. Xu, L. Sheng}, \textit{Positive semi-definite matrices of adjointable operators on Hilbert $C^*$-modules},
Linear Algebra Appl., \textbf{428}(2008), no. 4,   992--1000.

\bibitem{Xu-Wei-Gu}\textsc{Q. Xu, Y. Wei, Y. Gu}, \textit{Sharp norm estimations for Moore-Penrose inverses of stable
perturbations of Hilbert $C^*$-module operators}, SIAM J. Numer. Anal., \textbf{47}(2010), 4735--4758.

\bibitem{Xu-Yan}\textsc{Q. Xu, G. Yan}, \textit{Harmonious projections and Halmos' two projections theorem for Hilbert $C^*$-module operators}, Linear Algebra Appl., \textbf{601}(2020), 265--284.

\bibitem{XY2}\textsc{Q. Xu, G. Yan}, \textit{Products of projections, polar decompositions and norms of differences of two projections}, Bull. Iranian Math. Soc.,
\textbf{48}(2022), 279--293.







\end{thebibliography}
\end{document}